\documentclass[mathpazo]{cicp}

\usepackage{amsmath,epsf,cite}
\usepackage{amssymb,amsthm,tocvsec2}
\usepackage{bigints}
\usepackage{graphicx}
\usepackage{sidecap}
\usepackage{wrapfig}
\usepackage{pict2e}
\usepackage{float}
\usepackage{graphicx}
\usepackage{setspace}
\usepackage{float}
\usepackage{color}
\usepackage{mathrsfs}
\usepackage{epsfig,epsf,latexsym,subfigure}
\usepackage{fancybox}
\usepackage{fancyhdr}
\usepackage{mdwlist}
\usepackage{paralist}
\usepackage{setspace}
\usepackage{tikz}
\usetikzlibrary{shapes,arrows}
\usepackage{enumitem}
\setlist{nolistsep}
\usepackage{indentfirst}
\usepackage{epstopdf}

\renewcommand\arraystretch{1.5}

\numberwithin{equation}{section}
\numberwithin{figure}{section}
\numberwithin{table}{section}

\theoremstyle{plain}

\theoremstyle{remark}

\newcommand{\ben}{\begin{eqnarray}}
\newcommand{\een}{\end{eqnarray}}
\newcommand{\benx}{\begin{eqnarray*}}
\newcommand{\eenx}{\end{eqnarray*}}
\newcommand{\beq}{\begin{equation}}
\newcommand{\eeq}{\end{equation}}
\newcommand{\beqx}{\begin{equation*}}
\newcommand{\eeqx}{\end{equation*}}
\newcommand{\bea}{\begin{array}}
\newcommand{\eea}{\end{array}}
\newcommand{\bef}{\begin{figure}[H]}
\newcommand{\eef}{\end{figure}}
\newcommand{\be}{\begin{equation}}
\newcommand{\ee}{\end{equation}}

\newcommand{\bse}{\begin{subequations}}
\newcommand{\ese}{\end{subequations}}

\def\hR{\mathbb R}

\def\cF{\mathcal{F}}

\def\be{\mathbf{e}}

\def\bx{\mathbf{x}}

\def\cU{\mathcal{U}}
\def\cF{\mathcal{F}}

\begin{document}
\title{Solving Allen-Cahn and Cahn-Hilliard Equations using the Adaptive Physics Informed Neural Networks}

\author[Wight and Zhao]{ Colby L. Wight \affil{1} and
Jia Zhao\affil{1}\comma \corrauth}
\address{ \affilnum{1}\ Department of Mathematics \& Statistics, Utah State University, Logan, UT, USA }
\emails{{\tt colbywight5@gmail.com.} (C.~Wight), {\tt jia.zhao@usu.edu.} (J.~Zhao)}

\begin{abstract}
	Phase field models, in particular, the Allen-Cahn type and Cahn-Hilliard type equations, have been widely used to investigate interfacial dynamic problems. Designing accurate, efficient, and stable numerical algorithms for solving the phase field models has been an active field for decades. In this paper, we focus on using the deep neural network to design an automatic numerical solver for the Allen-Cahn and Cahn-Hilliard equations by proposing an improved physics informed neural network (PINN). Though the PINN has been embraced to investigate many differential equation problems, we find a direct application of the PINN in solving phase-field equations won't provide accurate solutions in many cases. Thus, we propose various techniques that add to the approximation power of the PINN. As a major contribution of this paper, we propose to embrace the adaptive idea in both space and time and introduce various sampling strategies, such that we are able to improve the efficiency and accuracy of the PINN on solving phase field equations. In addition, the improved PINN has no restriction on the explicit form of the PDEs, making it applicable to a wider class of PDE problems, and shedding light on numerical approximations of other PDEs in general.
\end{abstract}

\maketitle

\section{Introduction}
Phase field models have been widely embraced in the past few decades to study various problems in science and engineering, taking the applications in image analysis, material science, engineering, fluid mechanics, and life science as examples. Among them, two fundamental equations are the Allen-Cahn (AC) equation and Cahn-Hilliard (CH) equation, which are originally introduced to describe the non-conservative and conservative phase variables in the phase separation process, respectively. Both models are recognized as gradient flow systems, for which there exists a Lyapunov function, known as the free energy. From a modeling view, given a specified Lyapunov function, or a free energy function, the Allen-Cahn type equations can be derived as the $L^2$ gradient flow, and the Cahn-Hilliard type equations can be derived as the $H^{-1}$ gradient flow, respectively. This generality makes the AC and CH type equations extremely useful in modeling many interfacial or multiphase problems. And many well-known PDE models turn out to be their special cases.

Given the nonlinearity in phase field equations, along with the stiff terms due to the small parameters, how to design accurate, efficient, and stable numerical algorithms for their numerical approximations have been intensively studied in the literature. Here is some literature that attracts our attention \cite{Li&KimCNSNS2017, Guillen-JCM, ZQWdendri2017, Shen&Wang&Wise2012, HanWang2015, ShenJ2, Yang&ZhaoCICP2018}. Interested readers are encouraged to read them and the references therein for further information.
In this paper, we focus on a new numerical approximation approach by using the deep neural network. Our major goal is to investigate strategies to improve the capabilities of deep neural networks on solving phase field models, in particular, the Allen-Cahn equation and the Cahn-Hilliard equation.

The artificial neural network is named after the fundamental unit of computation inside the mammalian brain \cite{mehlig_2019}. Many neurons inside the brain work together to carry out complex tasks. Similarly, an artificial neural network is composed of multiple connected neurons that work to solve complex tasks. 
A single neuron in a neural network can take input from multiple neurons (or nodes). Each input has a parameter called a weight associated with it. There is also typically a bias term that doesn't have an input associated with it. The neuron receives the sum of these inputs multiplied by their weights, along with added bias. This weighted sum then goes through an activation function that gives the final output for this neuron. In the brain, a neuron usually doesn't fire unless the total of its input reaches a certain threshold. The output is either on or off. In deep learning continuous activation functions are more commonly used \cite{activation_fn}. The $sigmoid$ function can be used as a smoother version of the step function. There are benefits in using differentiable functions like this that will help in "learning" good weights. Other useful activation functions used in deep learning include $relu$, $\tanh$, leaky $relu$ \cite{leaky_relu} and $swish$ \cite{swish}.
A typical feed-forward, fully connected neural network has input going to and from multiple neurons. The input to the network makes up the input layer. The neurons of the input layer are then sent to other layers of neuron connections called hidden layers, and finally to the last layer, the output layer. See Figure \ref{fig:neural-network} for a representation of a simple neural network architecture.

\def\layersep{2.5cm}
\begin{figure}
	\centering
	\begin{tikzpicture}[shorten >=1pt,->,draw=black!50, node distance=\layersep]
	\tikzstyle{every pin edge}=[<-,shorten <=1pt]
	\tikzstyle{neuron}=[circle,fill=black!25,minimum size=17pt,inner sep=0pt]
	\tikzstyle{input neuron}=[neuron, fill=green!50];
	\tikzstyle{output neuron}=[neuron, fill=red!50];
	\tikzstyle{hidden neuron}=[neuron, fill=blue!50];
	\tikzstyle{annot} = [text width=4em, text centered]
	
	\foreach \name / \y in {1,...,4}
	\node[input neuron, pin=left:\#\y \,\,\, Input ] (I-\name) at (0,-\y) {};
	
	\foreach \name / \y in {1,...,5}
	\path[yshift=0.5cm]
	node[hidden neuron] (H-\name) at (\layersep,-\y cm) {};
	
	\node[output neuron,pin={[pin edge={->}]right:Single Output}, right of=H-3] (O) {};
	
	\foreach \source in {1,...,4}
	\foreach \dest in {1,...,5}
	\path (I-\source) edge (H-\dest);
	
	\foreach \source in {1,...,5}
	\path (H-\source) edge (O);
	
	\node[annot,above of=H-1, node distance=1cm] (hl) {Hidden layer};
	\node[annot,left of=hl] {Input layer};
	\node[annot,right of=hl] {Output layer};
	\end{tikzpicture}
	\caption{A diagram of a neural network with an input layer, hidden layer, and output layer. Each input gets sent to each neuron in the hidden layer. The arrows between the neurons all have a weight associated with them. The bias for each hidden neuron and output neuron are not shown. }
	\label{fig:neural-network}
\end{figure}
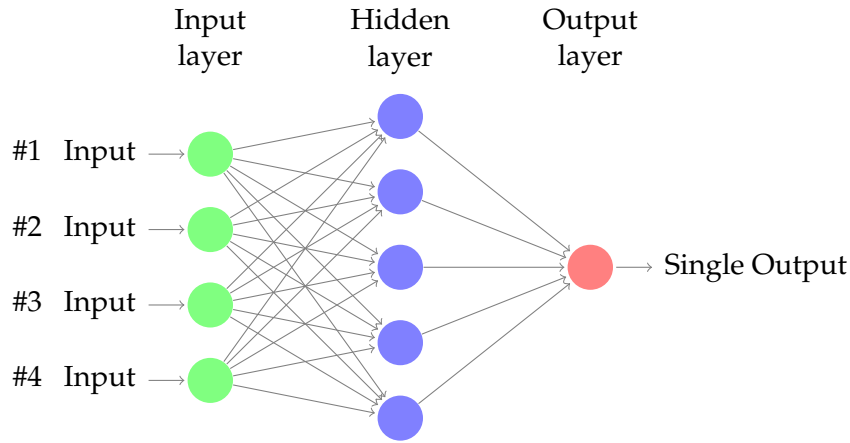

Mathematically, the feed-forward neural network could be defined as compositions of nonlinear functions. Given an input $x \in \hR^{n_1}$, and denote the output of the $l$-th layer as $a^{[l]} \in \hR^{n_l}$, which is the input for $l+1$-th layer. In general, we can define the neural network as \cite{DeepLearningMath} 
\beq \label{eq:NN}
\bea{l}
a^{[1]} = x \in \hR^{n_1}, \\
a^{[k]} = \sigma \Big( W^{[k]} a^{[k-1]} + b^{[k]} \Big) \in \hR^{n_k}, \quad \mbox{ for } k=2,3,\cdots, L,
\eea 
\eeq 
where  $W^{[k]} \in \hR^{n_k \times n_{k-1}}$ and $b^{[k]} \in \hR^{n_k}$ denote the weights and biases at layer $k$ respectively, $\sigma$ denotes the activation function.
Essentially neural networks are non-linear mappings with many parameters. Due to a large number of parameters, they are referred to as a "black-box" \cite{black_box}. A loss function is defined to measure the difference between the output of the network with the known (expected) output. Then, the parameters (the weights and biases) can be "learned" by minimizing the loss function. 
Typically a gradient-based optimization scheme is applied. The backpropagation algorithm is an efficient way to find the gradient of this highly dimensional loss function \cite{hecht1992theory}. Stochastic gradient descent is a popular method that uses smaller subsets of the training data called batches for each step to achieve better results \cite{bottou2010large}. 

In the past few years, there has been intensive research on understanding how deep neural networks can be adopted to solve and discover differential equations \cite{Adaptive,RaissiJCP,DeepXDE,DGM,Rudy-SIAM,PDE-discovery,multistep,DHPM,discussion-pde,complex-geometries,E,Yu}.
For instance, various studies deal with problems of data-driven modeling, especially ones using deep learning to solve differential equations \cite{HanEPNAS2018}. Others propose approaches to use existing experiment data to identify the differential equations themselves in order to model real-world phenomenon \cite{Brunton3932, rudy2016}.
In particular, Raissi et al. propose the physics informed neural network (PINN) to aid in both the solution of differential equations as well as their discovery \cite{RaissiJCP}. The PINN was shown to solve Burgers' equation and the Schrodinger equation with certain initial conditions accurately. Though the PINN has been widely appreciated in the community, we found it is not capable of solving the phase field equations in many cases properly.
As a result, we observe that some parts of solutions (both spatially and temporally) for the phase field models might be harder to learn than other parts. 
These difficult areas could even change over the course of learning the solution. This is reasonable given the solutions for phase field models usually have sharp transition layers, and the layers evolve with time. Unfortunately, the original PINN method \cite{RaissiJCP} failed to address these issues, as it involves choosing training data points at the beginning of training randomly across the domain and fixes them over the entire training process. 

This motivates us to conduct the research in this paper. This major goal of this paper is to investigate strategies to improve the approximation capabilities of physics informed deep neural networks on solving the phase field equations. We propose strategies to improve the accuracy and efficiency of PINNs. The major contribution of this paper is to involve the idea of adaptive sampling data points over the course of training, both with time and space position in mind. In the rest of this paper, we will introduce our ideas in detail, and provide several numerical examples to justify our approach further.

\section{Improved Physics Informed Neural Networks}
In this section, we will first give a brief review of the physics informed neural network (PINN) \cite{RaissiJCP}. Then, we propose several strategies to improve the accuracy and approximating capability of the PINN.

\subsection{Physics informed neural networks}

To elucidate the idea of the physics informed neural network (PINN) \cite{RaissiJCP}, we use the following Burgers' equations as an example, shown as
\begin{equation}
\begin{split}
&u_t+uu_x-(0.01/\pi)u_{xx}=0, \quad x\in[-1,1],\quad t\in[0,1], \\
&u(0,x)=-sin(\pi x), \\
&u(t,-1)=u(t,1)=0.
\end{split}
\end{equation}
To solve the Burgers' equation, the authors \cite{RaissiJCP} introduce two neural networks: the $u$-network as 
\beq
\mathcal{U}: (x, t) \rightarrow  \mathcal{U}(x,t),
\eeq 
and the $f$-network as
\beq
\mathcal{F}: (x, t) \rightarrow  \cU_t(x,t) + \cU(x,t)  \cU_x(x,t)  - (0.01/\pi)\cU_{xx}(x,t) .
\eeq 
In this case, both the $u$ and $f$ networks have two input neurons, and one output in the final (output) layer.
Given the initial data points $\left\{ (0, x_u^i, u^i) \right\}_{i=1}^{N_u}$ and random samples from the boundary $\left\{ (t_b^i, -1), (t_b^i, 1)\right\}_{i=1}^{N_b}$ and collocation points  $\left\{ (t_f^i, x_f^i) \right\}_{i=1}^{N_f}$ in the interior of the domain $[-1, 1]\times[0, 1]$, the loss function for this problem is defined as the mean squared error of the $u$-network plus the mean squared error of the $f$-network, i.e.
\beq  \label{eq:Burgers-loss}
MSE = MSE_u+MSE_b + MSE_f,
\eeq 
where the three terms are defined as
\beq \label{eq:Burgers-loss-terms}
\begin{split}
	&MSE_u=\frac{1}{N_u}\sum_{i=1}^{N_u}|\cU(0,x^i_u)-u^i|^2, \\
	& MSE_b=\frac{1}{N_b}\sum_{i=1}^{N_b}|\cU(t^i_u,1)|^2 + |\cU(t^i_u,-1)|^2, \\
	&MSE_f=\frac{1}{N_f}\sum_{i=1}^{N_f}|\cF(t^i_f,x^i_f)|^2,
\end{split}
\eeq 
where $N_u$, $N_b$ and $N_f$ are the number of initial training data, boundary training data, and interior collocation points respectively. 
Here $(0, x_u^i)$ are the initial condition points that serve as the inputs to the u-network and $u^i$ is the actual solution value of $u$ at those points (provided in the initial condition). $t_f^i$ and $x_f^i$ are the collocation points passed into the $f$-network. Given the loss function value in \eqref{eq:Burgers-loss} is small enough, the $u$-network approximates the solution well. 

\subsection{Some strategies to improve the accuracy of the physics informed neural networks}
Though the PINN is shown to be powerful to solve the Burgers' equation, we find that a direct application of the PINN on solving the phase field equations would not provide accurate solutions in many cases. This motivates us to seek some extra techniques to improve the accuracy of the PINN.  In this section, we introduce some strategies to improve the approximation power of the PINN.

\subsubsection{Adding weights in the loss function} \label{sec:weights}
One simple technique to improve the learning capability of the PINN for solving phase field models is to add weights in the loss function. This is motivated by the fact that the phase field equations are dissipative, i.e., not reversible. For example, the Allen-Cahn equation, as with other reactive diffusion equations, can only be solved in the forward time direction. In other words, if the PINN fails to learn the solution at time $t=t_1$ well, there is little hope of learning the solution at a later time $t=t_2$ (with $t_2>t_1$) accurately. In order to put an emphasis on the importance of first learning the solution near $t=0$, we put more weight on $MSE_u$ to enforce the $u$-neural network satisfies the initial condition. Specifically, we can redefine the loss function in \eqref{eq:Burgers-loss} as follows:
\beq \label{eq:loss-weight}
MSE = C_0 MSE_u+MSE_b+MSE_f,
\eeq 
where $C_0$ is a big positive constant (that is tunable). Here we chose $C_0=100$ is this paper, if not specified separately.

\subsubsection{Mini-batching strategy to improve convergence} \label{sec:mini-batch}
Mini-batching is a technique that has been used in deep learning to improve performance. Instead of using the entire data set to calculate the exact direction of the gradient, a subset of the data called a batch or mini-batch is used to evaluate the direction. Mini-batching has been shown to help avoid less desirable local minimum better than full-batch gradient descent \cite{kleinberg2018}. 

We point out that, in \cite{RaissiJCP} and many follow-up papers, the authors do not use a mini-batching strategy in their training process. In this paper, we investigate the mini-batching approach and observe the mini-batch approach indeed can facilitate the convergence of the trained neural network for approximating phase field equations.

\subsection{Adaptive strategies to improve the accuracy of the physics informed neural networks}
In this section, we borrow the idea of temporal and spatial adaptivity in classical numerical methods for solving differential equations, and introduce some adaptive strategies to improve the accuracy of the physics informed neural networks.

\subsubsection{Adaptive sampling of collocation points} \label{sec:adaptive-sample}
Instead of using fixed sample points over the training process, we realize an adaptive resampling of the collocation points across the domain during the training process is essential for certain situations. In particular, for the phase field equations, there are moving interfaces that are sharp transitions over space, where finer meshes are desired to capture the evolution dynamics.
Therefore, instead of only sampling points evenly across the domain, we periodically stop training and re-evaluate where points are needed most. 
We notice there is a correlation between the points that had a larger error in the $u$-network (the solution) and the points that had large errors in the $f$-network. This motivates us to use the error of $f$-network as an indicator for resampling.

\begin{figure}
\centering
\includegraphics[width=\textwidth]{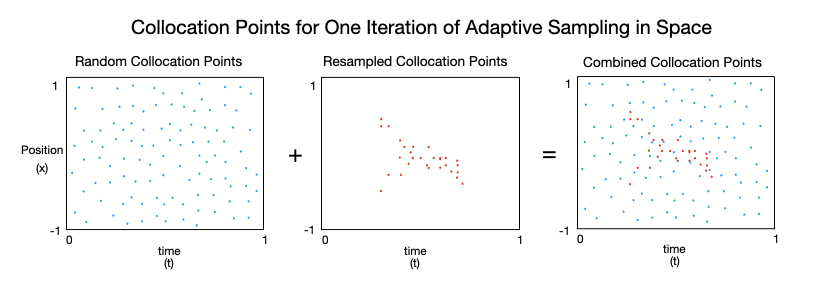}
\caption{An illustration of adaptive sampling the collocation points over one iteration. Consider the domain $x \in [-1, 1]$ and $t \in [0, 1]$. The blue points show the set of randomly sampled collocation points using Latin hypercube sampling. Training on these points keeps the solution of the equation accurate across the whole domain. The red points show an example of a set of resampled collocation points sampled after evaluating the $f$-prediction network for the highest areas of error. These points help to improve the solution accuracy over red-point zones, which is usually the interfaces for phase field equations. The combined set of points is the collocation points used to train the network. The network can repeat this process for multiple re-sampling iterations. The blue points will stay the same, but the red points may change to focus on other parts of the domain that are not being learned well. }
\label{fig:sapce-sample}
\end{figure}

In practice, we first train the network using the randomly selected points across the domain. We then choose a different set of sample test points across the domain using the same Latin hypercube sampling technique and pick a portion of the points that give the highest error in predicting $f$. We add this set of points to previous collocation points and train the network again. This procedure is important, as it prevents us from losing the accuracy of the solution across the entire domain and, in the meanwhile, helps us to focus more points to learn the trickier parts better. This process can be iterated a few times if necessary, by adding a different set of sampled collocation points to the original set and training again. An illustrative diagram is shown in Figure \ref{fig:sapce-sample}.

\subsubsection{Adaptive strategies in time} \label{sec:adaptive-in-time}
The strategy of adaptive sampling collocation points addresses the moving interface for the phase field solutions. But for certain phase field problems with sharp transition, even the PINN with an adaptive sampling of collocation points fails to converge to the actual solution. Thus, some extra attention is still needed.

In this section, we introduce two adaptive strategies in time to improve the convergence of the PINN.  The first time-adaptive approach is similar to the adaptive method introduced above, in that collocation points in time are strategically chosen to improve learning. The second time-adaptive method takes a different approach where we create separate networks on smaller (subsequent) time domains of fixed or adaptive length.  

\begin{figure}
\centering
\includegraphics[width=\textwidth]{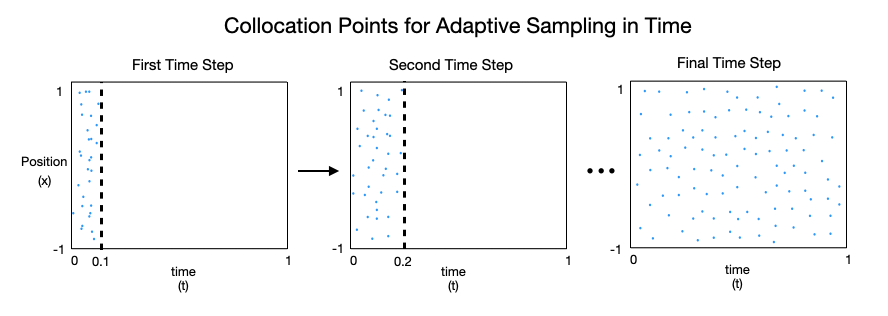}
\caption{ This figure illustrates the concept behind the time-adaptive approach I. Consider a time domain $[0, 1]$. The first time step only allows data points to be taken between $t \in [0, 0.1]$. Once the $f$-predicted error on the interval is sufficiently small, collocation points are then chosen on a larger interval $[0, 0.2]$, including adaptive space sampling. Collocation points are still chosen from the earlier time domains to keep what has been learned there learned well. This is continued until the time interval covers the entire time range for the problem. Note this is all done on one PINN (in the next time method, multiple networks are created for each time interval). }
\label{fig:time-adaptive-approach-I}
\end{figure}

\textbf{Time-adaptive approach I: adaptive sampling in time.}
At each time step of this approach, we require the data points, (initial, boundary, and general collocation both original and resampled) to come from within a specified time interval. For instance, if we are approximating the solution in the time domain $[0, 1]$, we start with small time intervals $[0, t_1]$, $t_1 >0$, where $t_1$ is close to zero, saying $t_1=0.1$. Then we gradually increase the time span, i.e., $[0, t_i]$, $i=1,2,\cdots, N$, with $0<t_1<t_2< \cdots < t_N = 1$, when each time span is learned well. Eventually, the solution is learned well on the whole domain. This idea is illustrated in Figure \ref{fig:time-adaptive-approach-I}. In particular, one usually can set a threshold and a maximum training iterations. Once the loss function value is smaller than the threshold or the training exceeds the maximum training iterations, the training process marches onto next time step. In certain cases, if the loss function value is still huge after maximum training iterations, the time step size shall be reduced.

In addition, this idea of adaptive sampling in time allows the user to designate a list of time steps. For each time step, collocation points will be sampled from only this restricted domain. Adaptive space sampling is used within the restricted time interval to improve learning. The network is then trained for each time interval using adaptive sampling collocation points. At each new iteration for the time step, the $f$-predicted error on that time domain is calculated.

\textbf{Time-adaptive approach II: adaptive time marching strategy.}
In the second time-adaptive approach, we propose
to split up the domain of interest into smaller problems. Notice that in the first time-adaptive approach, we only have a single network that focuses the collocation points adaptively in time. In the second approach, we create separate networks for each time step (interval). For example, if our domain of interest is $[0, 1]$, and we set the time step $\Delta t = 0.1$. we train one network to learn the solution on the interval from $[0, 0.1]$. Once the solution is learned well on this time interval, we train another network on the interval $[0.1, 0.2]$, and so on, until we solve the problem in the entire time domain.

\begin{figure}
\centering
\includegraphics[width=\textwidth]{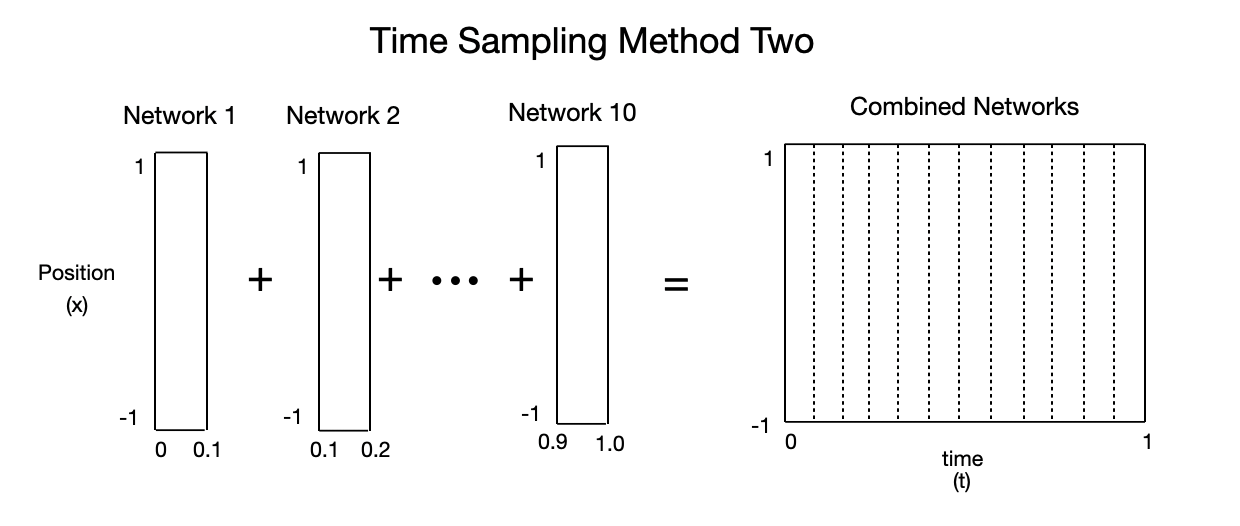}
\caption{ This figure illustrates the concept behind the adaptive time-marching strategy. Here we choose time step $\Delta t =0.1$, and individual networks are trained for each time step. Network 2 shares the same initial time as the last time in Network 1. Once Network 1 has been learned well, its values at $t=0.1$ can be used as the initial condition for Network 2. The solution for each individual network can be combined at the end into one continuous solution that covers the entire domain of the problem. Each network has the same time length but can be bigger or smaller for parts of the domain that are easier or harder to learn. }
\label{fig:time-adaptive-approach-II}
\end{figure}

A caveat here is we cannot use the initial condition for the later networks as the given initial value is only valid for $t=0$. Instead, we use the prediction from the previous time step's network as the initial condition for the current network. We continue doing this until we have covered the whole time domain of the original problem. The individual networks can be combined to obtain a solution at any point in the domain of the original problem. The idea is illustrated in Figure \ref{fig:time-adaptive-approach-II}.

We remark that, for the solutions of phase field models, there are sharp transitions in both time and space, which make them very tough to solve. Utilizing the adaptive ideas (in both time and space),  our improved (adaptive) PINNs improve the accuracy of the baseline PINN by avoiding certain local minima or saddle points.

\section{Numerical Results}
In this section, we will provide several numerical tests on solving the Allen-Cahn equation and the Cahn-Hilliard equation with the proposed strategies in the previous section.

First of all, let us recall the Allen-Cahn equation and Cahn-Hilliard equation. Denote the specific expression of the Ginzburg-Landau free energy
\beq
F = \int_\omega \frac{\gamma_1}{2}|\nabla u|^2 + \frac{\gamma_2}{4}(u^2-1) d\bx,
\eeq 
where $\gamma_1$ and $\gamma_2$ are parameters. If we take the $L^2$ gradient flow, we obtain the Allen-Cahn equation
\beq \label{eq:AC}
\partial_t u = \gamma_1 \Delta u + \gamma_2 (u - u^3);
\eeq 
and if we take the $H^{-1}$ gradient flow, we obtain the Cahn-Hilliard equation
\beq \label{eq:CH}
\partial_t u = \Delta ( - \gamma_1 \Delta u + \gamma_2 (u^3-u)).
\eeq 
Thermodynamically consistent boundary conditions (such as periodic boundary conditions and homogeneous Neumann boundary conditions) and initial values shall be proposed to close the system.



In the rest of this section, we use uniformly random samples from the initial and boundary points and use a Latin hypercube sampling (LHS) strategy to sample the collocation points in the interior domain. 
Recall the $u$-network gives the approximated value of the solution to the differential equation, and the $f$-network gives the approximated residual value of the equation at the given point, which should ideally be zero.

To test the accuracy of the learned solution, the \lq{}actual\rq{} solution is obtained by solving the phase field equations using classical numerical methods with high accuracy. The accuracy of the trained model is assessed by taking the relative $l_2$-norm 
of the difference between the \lq{}actual\rq{} value $u(x_i, t_i)$ at those points and the $u$-network output $\cU(x_i, t_i)$ at those points, i.e.,
\beq
Error = \frac{\sqrt{\sum_{i=1}^N | \cU(x_i,t_i) - u(x_i, t_i)|^2}}{ \sqrt{\sum_{i=1}^N| u(x_i, t_i)|^2}},
\eeq 
given the data points $\{(x_i, t_i) \}_{i=1}^N$, with $N$ the number of the points.

\subsection{Solving the Allen-Cahn equation}
We first tested the Allen-Cahn equation of \eqref{eq:AC} with periodic boundary conditions in one dimension, and chose the parameters: $\gamma_1=0.0001$ and $\gamma_2 = 5$. The specific system is summarizsed as follows:
\begin{equation} \label{eq:AC-test-1}
\begin{split}
&u_t - 0.0001 u_{xx}+5u^3-5u=0,\quad x \in [-1,1],\quad t\in [0,1], \\
&u(0,x) = x^2cos(\pi x), \\
&u(t,-1)=u(t,1), \\
&u_x(t,-1)=u_x(t,1).
\end{split}
\end{equation}

Note that in \cite{RaissiJCP}, the authors did not test this problem with the (continuous) PINN. Instead, they solved it using the discrete Runge-Kutta neural network. We thus first test it with the PINN from \cite{RaissiJCP} as a baseline. 

As mentioned in previous section, we introduce two neural networks: the $u$-network as
\beq
\mathcal{U}: (x, t) \rightarrow  \mathcal{U}(x,t),
\eeq 
and the $f$-network as
\beq
\mathcal{F}: (x, t) \rightarrow  \cU_t(x,t)  - 0.0001 \cU_{xx}(x,t) + 5 \cU^3(x,t) - 5 \cU(x,t) .
\eeq 
And the loss function is defined as
\beq 
MSE = MSE_u+MSE_f+MSE_b,
\eeq 
Where $MSE_u$ and $MSE_f$ are defined the same as \eqref{eq:Burgers-loss-terms} and the error due to boundary term is replaced by 
\beq 
MSE_b=\frac{1}{N_b}\sum_{i=1}^{N_b} \Big( |\cU(t_b^i,x_u)-\cU(t_b^i,x_l)|^2+|\cU_x(t_b^i,x_u)-\cU_x(t_b^i,x_l)|^2 \Big)
\eeq
to address the periodic boundary condition.
Here $N_b$ is the number of collocation points used on the boundary, $\{t_b^i\}_{i=1}^{N_b}$ are the time values for those points, $x_u$ is the upper bound for $x$, and $x_l$ is the lower bound. For this problem we have $x_u = 1$ and $x_l = -1$ from our domain. We see that both the expressions inside the absolute values are ideally zero if the learned solution follows the boundary conditions set by the equation. In this problem (otherwise specified), we use $N_f=10,000$ collocation points, $N_u=200$ initial points, $N_b=200$ boundary points,
and $\tanh$ as the activation function. 
As a general practice in the rest of this paper, to optimize the loss function, we use Adam optimizer first and then use the L-BFGS-B optimizer to fine-tune the neural network.

\begin{figure}
\centering
\includegraphics[scale=.65]{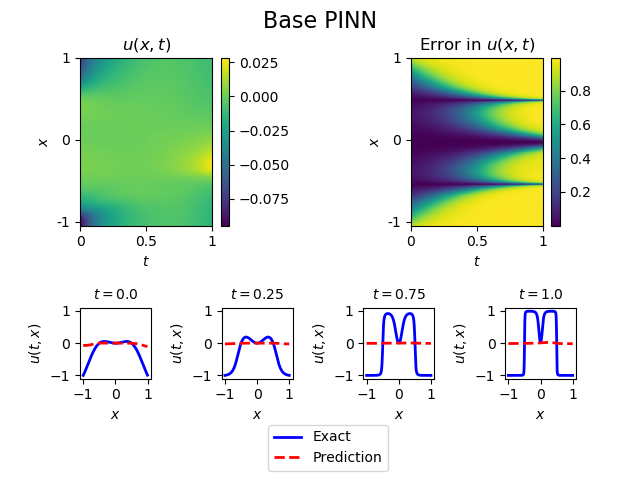}
\caption{Solutions of the AC equation learned using the Base PINN method. The four plots on the bottom are the predicted solutions vs. accurate solutions at different times. The baseline PINN approach failed to learn the solution well. }
\label{fig:PINN-AC}
\end{figure}

First of all, we tested the baseline PINN approach \cite{RaissiJCP}.
Unfortunately, using the baseline PINN approach alone, we were not able to learn the accurate solution for the Allen-Cahn equation.  The relative $l_2$ error stayed around $0.99$. The result is summarized in Figure \ref{fig:PINN-AC}, and we observed that the predicted solution at different time steps is not close to the actual solution. 

\begin{figure}
\centering
\includegraphics[scale=.65]{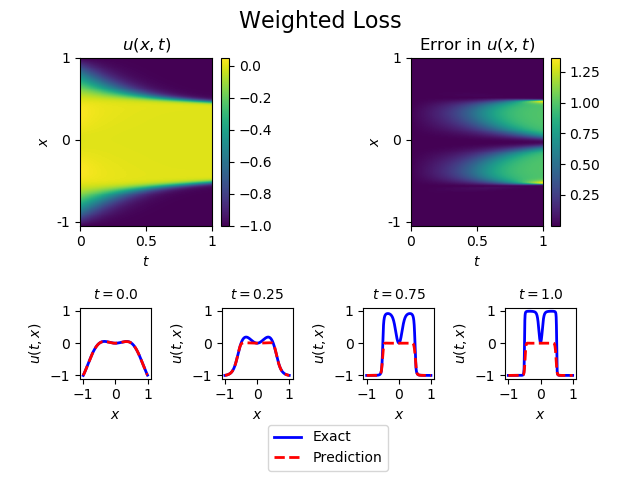}
\caption{Solutions of the AC equation learned using the PINN method with weights in the loss function. Here we were still using the baseline PINN approach without changing the collocation point, but we added weights in the loss function by putting more weight on the initial condition vs. the collocation and boundary conditions. This results in a slight improvement, especially at times near $t=0$. The learned solution is still not accurate enough.}
\label{fig:PINN-AC-weight}
\end{figure}

Then, we added the weights in the loss function, as proposed in Eq.  \eqref{eq:loss-weight}. With this strategy, the neural network performed slightly better with a relative $l_2$ error of $0.52$, but the algorithm still failed to converge, as shown in Figure \ref{fig:PINN-AC-weight}. The network learned the solution better at points at the starting time domain and near the boundary, but the learned solutions lost their accuracy in the later time domain, where the collocation points are sampled randomly.

\begin{figure}
\centering
\includegraphics[scale=.65]{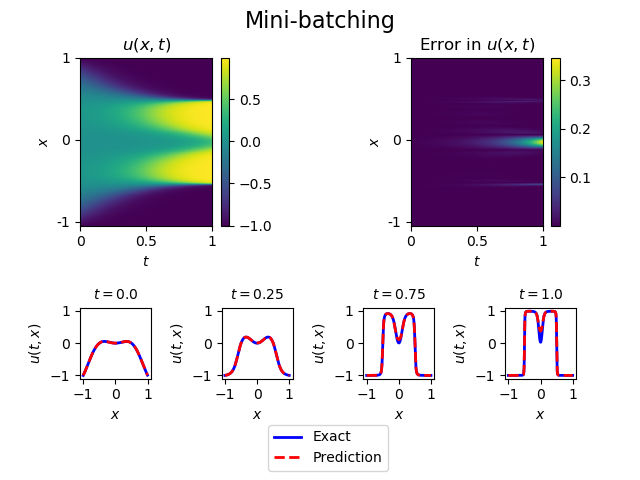}
\caption{Solutions of the AC equation learned using the PINN method being trained with mini-batching approach. The mini-batching approach helped improve the accuracy more than adding weights in the loss function alone. }
\label{fig:mini-batch-AC}
\end{figure}

Then, we used the mini-batching strategy introduced in Section \ref{sec:mini-batch}. In this trial, we used 10,000 collocation points, 512 initial points, and 200 boundary points.  This network had 4 hidden layers with 128 neurons per layer, and $\tanh$ as the activation function. 
The results are summarized in Figure \ref{fig:mini-batch-AC}, where we chose mini-batch size as $32$ with $100$ epochs.
The mini-batching trial performed almost as well as the sampling strategy in the previous trial. Notice, however, for both these solutions, the learned solutions near $t=1$ and $x=0$ were not exactly matching the actual solutions. By using time sampling approaches introduced in \ref{sec:adaptive-in-time}, the solutions could be learned accurately near both $t=1$ and $x=0$. 

\begin{figure}
\centering
\includegraphics[scale=.65]{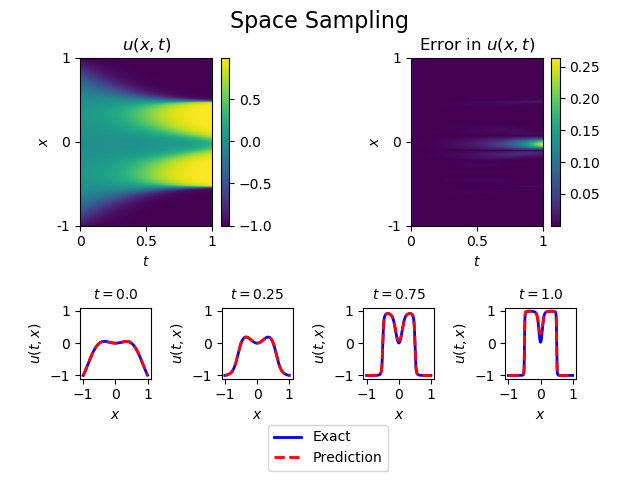}
\caption{Solutions of the AC equation learned using the PINN method with weighted loss function and adaptive sampling collocation points during the iteration. This improved PINN approach learned an accurate solution for the AC equation. The focusing of sampling collocation points in the areas with less accuracy allowed the network to learn the solution better. }
\label{fig:adaptive-sampling}
\end{figure}

Next, we used the adaptive sampling approach introduced in Section \ref{sec:adaptive-sample}. The network architecture was the same as the previous two trials. In this trial, we used 2000 original collocation points, along with 200 resampled collocation points for each iteration.
With adaptive sampling collocation points, we obtained a much better solution with higher accuracy. 
After about six re-sampling iterations, the solution did not improve significantly. The final result is shown in Figure \ref{fig:adaptive-sampling}. While the solution was much better, at the latter time steps, the learned solution didn't quite match the real solution. The relative $l^2$ error was improved from the others at $2.33e-02$. 

The $l_2$ errors for various cases are summarized in Table  \ref{tab:AC-error}. It is important to highlight also that this improved accuracy for the adaptive approach was obtained with only a fraction of the collocation points that the non-adaptive tests used. The adaptive test used only 2,000 collocation points while the others were tested with typical 10,000 collocation points. Thus, the adaptive approach noticeably reduces computational costs. We remark that the idea of adaptive re-sampling points is essential, as this same iteration scheme was tried without adding resampled points, and the solution did not converge as expected. 

\begin{table} 
\centering
\def\arraystretch{1.3}
\begin{center}
\begin{tabular}{|c|c|c|c|c|}
\hline
Allen-Cahn & PINN & Weighted Loss & Mini-batching & Re-sampling    \\
\hline
Relative $l_2$ & 9.90e-1 & 5.22e-1 & 3.25e-2 & 2.33e-2 \\ 
\hline
Relative $l_1$ & 9.90e-1 & 3.25e-1 & 8.80e-3 & 6.20e-3 \\ 
\hline
$l_\infty$-norm & 9.96e-1 & 1.37 & 3.37e-1 & 2.64e-1 \\
\hline
\end{tabular} 
\end{center}
\caption{Comparison of errors in the learned solutions of the Allen-Cahn equation using various PINN approaches: (1) the baseline  PINN approach; (2) adding weights in the loss function; (3) adoptive resampling collocation points. The adaptive re-sampling approach produced the best result. }
\label{tab:AC-error}
\end{table}

In the first example, we test the improved PINN method on solving the AC equation with a smooth interface. Next, we further investigate the capability of improved PINN on solving situations with sharper moving interfaces. In particular, we change the initial condition for the Allen-Cahn equation with more oscillations and use various values for the parameters to see how the methods work on different problems. 

In the following series of tests, we used the AC equation as below
\begin{equation}
\begin{split}
&u_t - \gamma_1 u_{xx}+\gamma_2u^3-\gamma_2u=0,\quad x \in [-1,1],\quad t\in [0,1], \\
&u(0,x) = x^2sin(2\pi x), \\
&u(t,-1)=u(t,1), \\
&u_x(t,-1)=u_x(t,1).
\end{split}
\end{equation}
This equation is different from the one previously tested of \eqref{eq:AC-test-1} in the initial condition and that the gamma parameters are not set. Instead of testing one set of parameters, we varied them, in particular $\gamma_2$, to see how our method performs on problems of increasing difficulty. We kept $\gamma_1=0.0001$ as in the previous problem. 

\begin{figure}
\centering
\includegraphics[scale=.65]{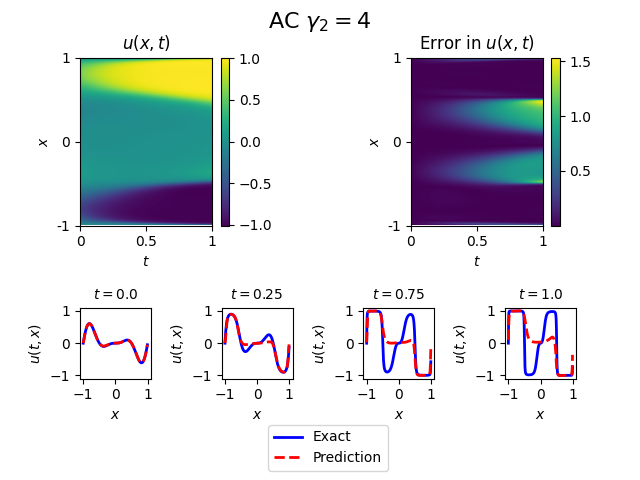} 
\caption{Solution results of the AC ($\gamma_2 = 4$). For this value of $\gamma_2$, we see that the solution is not learned all the way even with adaptive sampling.}
\label{fig:AC-gamma2-4}
\end{figure}

We tested the adaptive re-sampling method using parameters $\gamma_2 = 1, 2, 3, 4$. Note when $\gamma_2$ increases, the transition interface of the solutions is sharper, which makes it harder to solve the AC equation numerically. As an agreement, we observed numerically that for smaller values of $\gamma_2$, such as $\gamma_2=1$ and $\gamma_2=2$, the proposed resampling method in Section \ref{sec:adaptive-sample} converged in a reasonable number of re-sampling iterations. However, for $\gamma_2 = 3$, the resampling approach in Section \ref{sec:adaptive-sample} failed to converge within a reasonable number of re-sampling iterations. Note that this happened in fairly early time steps, before $t = .35$, and the error propagated and enlarged later on at time $t = 1$. It could be that a certain combination of the number of original collocation points with a number of resampled collocation points would help this solution to converge. However, we did not find one. This issue got even worse for $\gamma_2 = 4$. The neural network as unable to learn the larger curves in this case. 
As seen in Figure \ref{fig:AC-gamma2-4}, the sharp transition layer was not well captured by the learned solution.

\begin{figure}
\centering
\includegraphics[scale=.65]{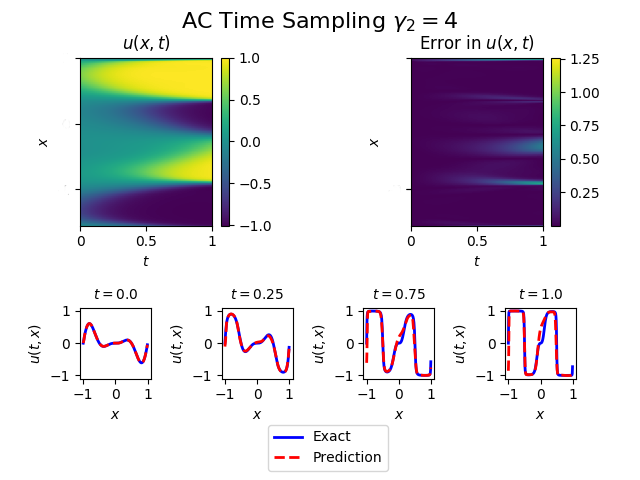}
\caption{Solutions of the AC equation ($\gamma_2=4$) solved using the improved PINN with the time-adaptive approach I and adaptive sampling of collocation points. This equation could not be solved using any of the methods used in previous trails. With fixed time steps of length $\Delta t =0.1$, this time-adaptive approach I focused on earlier times and then kept expanding the interval to encompass the whole domain. }
\label{fig:adaptive-I-AC}
\end{figure}

To overcome this, we utilized the adaptive in time sampling strategies introduced in Section \ref{sec:adaptive-in-time}. Finally, we saw the results for the trial using adaptive in time sampling approach I, as shown in Figure \ref{fig:adaptive-I-AC}. 
Using this approach, we were able to learn the solution for the more difficult problem with $\gamma_2 = 4$ much better with a relative $l_2$ error of $0.04$. As seen in Figure \ref{fig:adaptive-I-AC}, the learned solutions are accurate enough across the domain. As has been observed in other tests, the sharp curve when time is close to 1 is not learned perfectly. 

\begin{figure}[H]
\centering
\includegraphics[scale=.65]{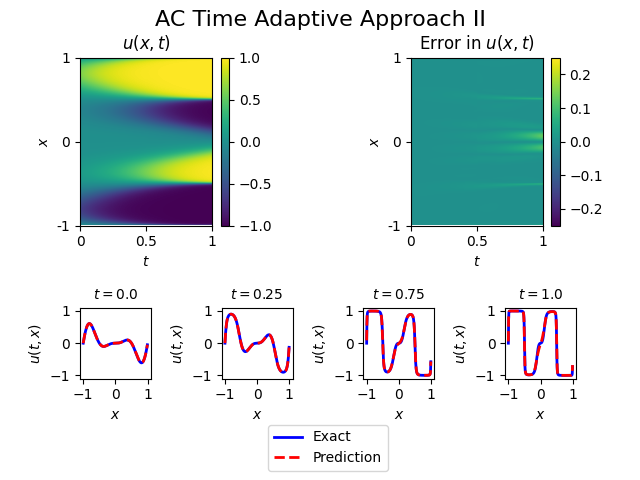}
\caption{Solutions of the AC equation ($\gamma_2=4$) solved using the improved PINN with the time-adaptive approach II and mini-batch training. With fixed time steps of length $\Delta t =0.25$, this time-adaptive approach II using multi-neural-networks to encompass the whole domain provides better accuracy. }
\label{fig:adaptive-II-AC}
\end{figure}

With time sampling approach II, we actually observed that difference to be smaller, as shown in Figure \ref{fig:adaptive-II-AC}. For the time sampling approach II on the Allen-Cahn equation, we saw that at times near $t=1$, the solution is learned better than it was in the past. 

Next, we used the proposed time-adaptive approach to further solve some more complicated problems in higher dimensions. We focused on the classical benchmark problem: the shrinking of a single drop problem for the Allen-Cahn equation.

\begin{figure}[H]
\centering
\subfigure[Predicted numerical solutions ($u_{pred}$)  at $t=0, 2.5, 5, 10$.]{
\includegraphics[width=0.25\textwidth]{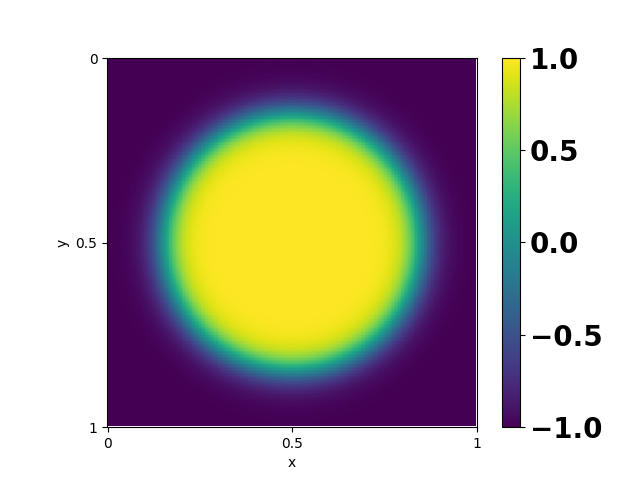}
\includegraphics[width=0.25\textwidth]{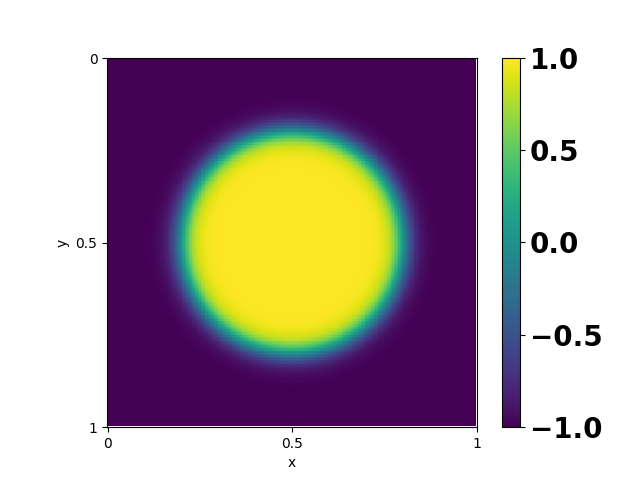}
\includegraphics[width=0.25\textwidth]{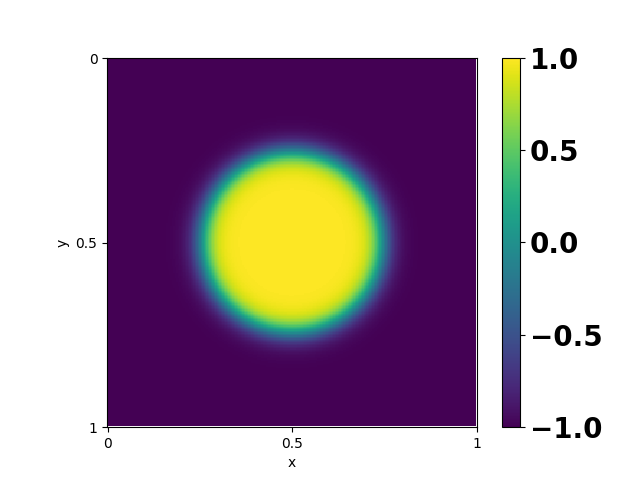}
\includegraphics[width=0.25\textwidth]{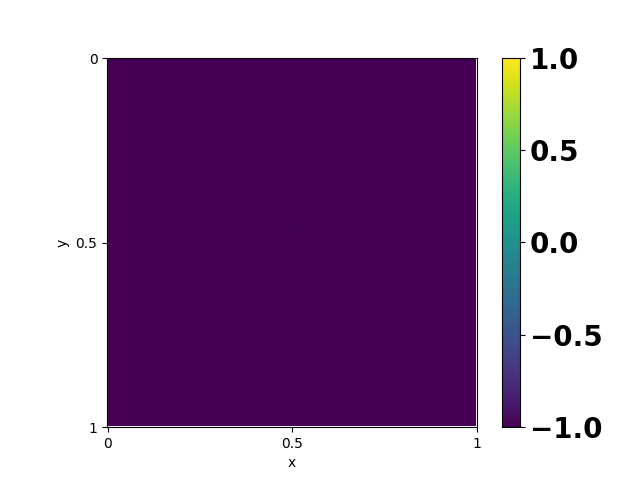}
}

\subfigure[Numerical errors ($u_{real}-u_{pred}$)  at  $t=0, 2.5, 5$ and $10$.]{
\includegraphics[width=0.25\textwidth]{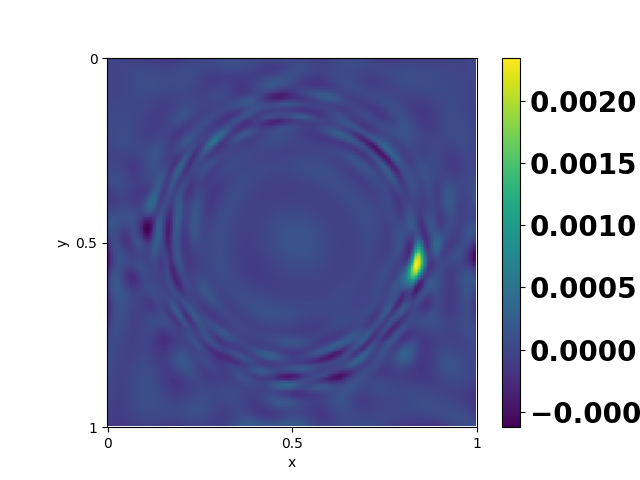}
\includegraphics[width=0.25\textwidth]{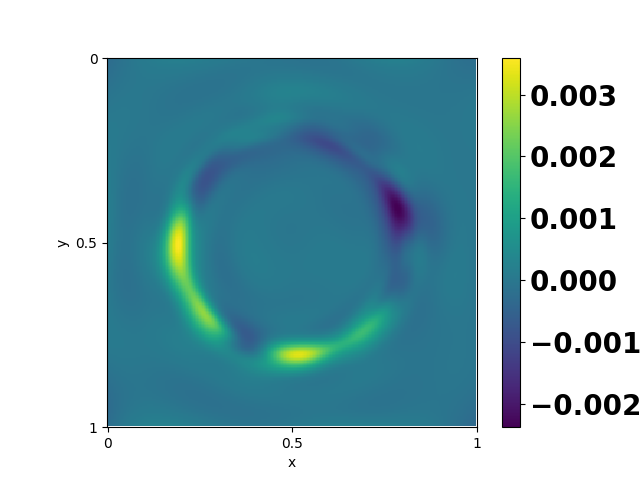}
\includegraphics[width=0.25\textwidth]{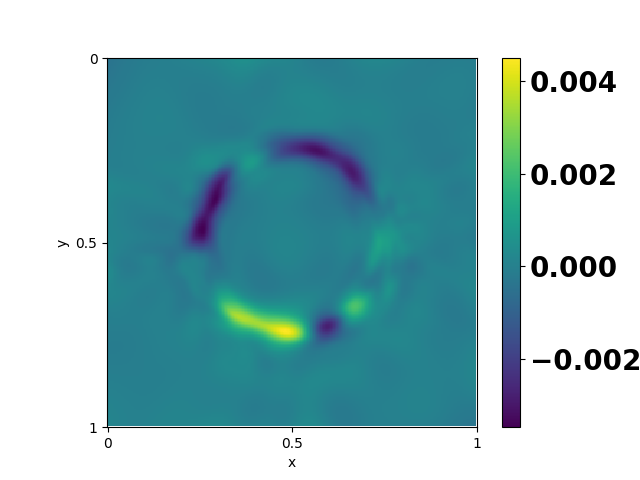}
\includegraphics[width=0.25\textwidth]{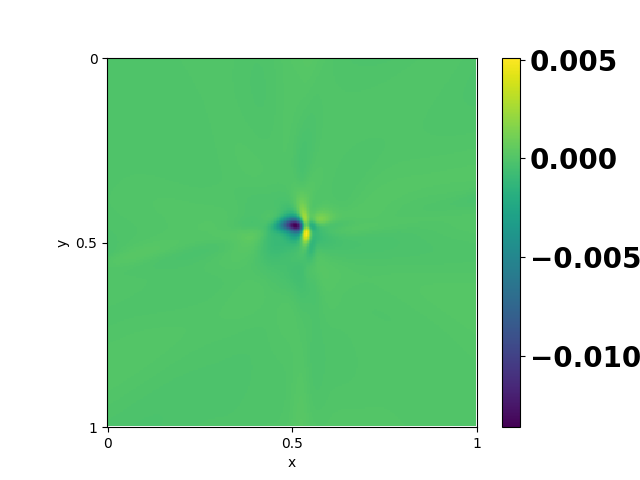}
}
\caption{Numerical approximations of 2D Allen-Cahn benchmark problem using the time-adaptive approach II. In this example, we chose $\Delta t=1$, and $C_0=10^3$   the neural network is $\cU: (x,y,t) \rightarrow \cU(x,y,t)$ with 6 hidden layers, 128 nodes per hidden layer. We chose time step $\Delta t = 0.2$. For each neural network, we used 100 epochs with batch size 32 for the Adam training, followed with a L-BFGS-B optimizer.}
\label{fig:AC2d}
\end{figure}

For the first benchmark problem, the Allen-Cahn model reads as
\beq
\partial_t \phi  = \lambda ( \varepsilon^2 \Delta \phi - \phi^3 +\phi), \quad \mathbf{X} \in \Omega, t\geq 0,
\eeq 
We chose  the domain $\Omega:=[0, 1]^2$, the parameters $\lambda=10$, $\varepsilon=0.025$. the initial profile for $\phi$ was given as
\beq
\phi(x, y, t= 0) = \tanh( \frac{ 0.35 - \sqrt{(x-0.5)^2 +(y-0.5)^2}}{2 \varepsilon}).
\eeq
and We utilized the time-adaptive approach II to solve this problem for $t \in [0, 10]$ with fixed time step size $\Delta t=1$. The predicted numerical solutions are shown in Figure \ref{fig:AC2d}(a), where we observed the drop shrink and eventually disappeared. The numerical error, i.e. the difference between the real solution (which is computed accurate with classical numerical solver) and the predicted solutions of the neural network are shown in Figure \ref{fig:AC2d}(b), where we observed the error was already very small. In other words, the time-adaptive approach II provided accurate approximation for this problem.

\begin{figure}[H]
\centering
\subfigure[Predicted numerical solutions ($u_{pred}$)  at $t=0, 0.2, 0.5, 1$]{
\includegraphics[width=0.25\textwidth]{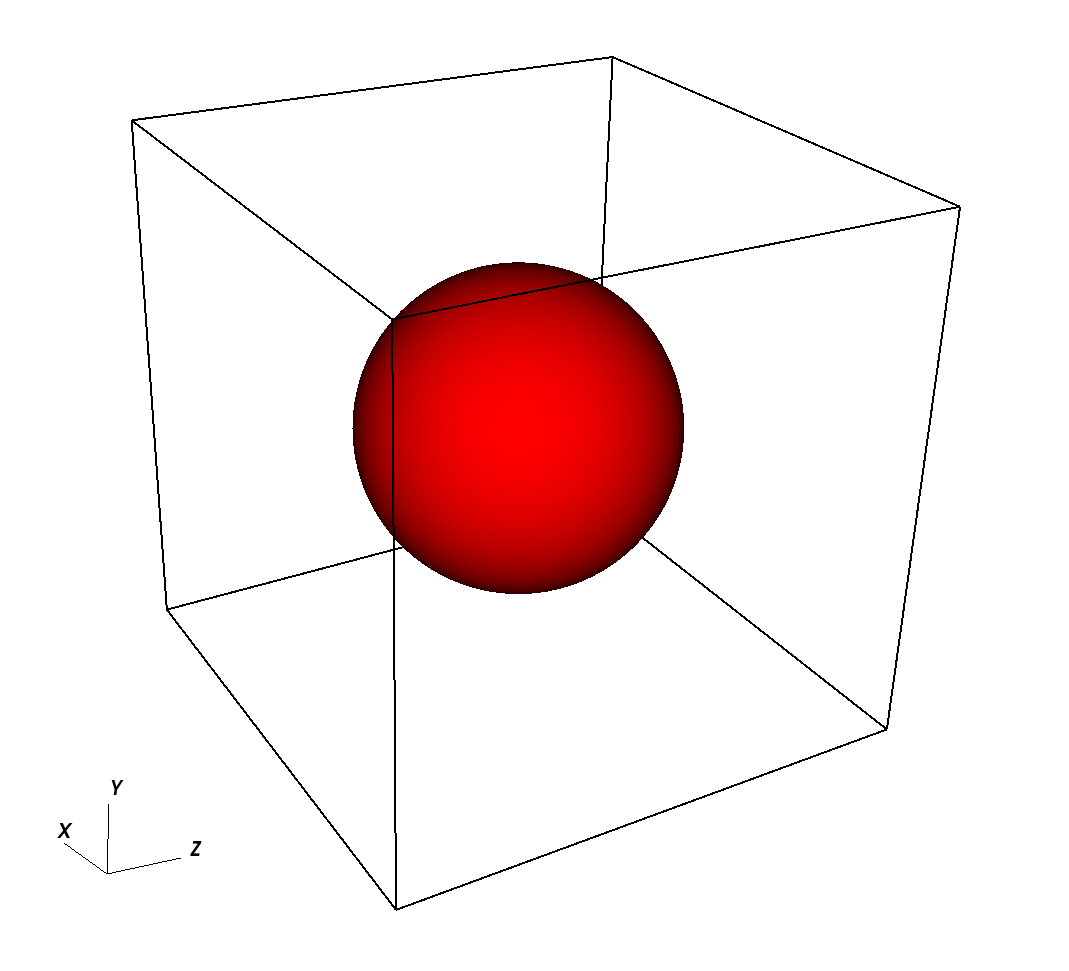}
\includegraphics[width=0.25\textwidth]{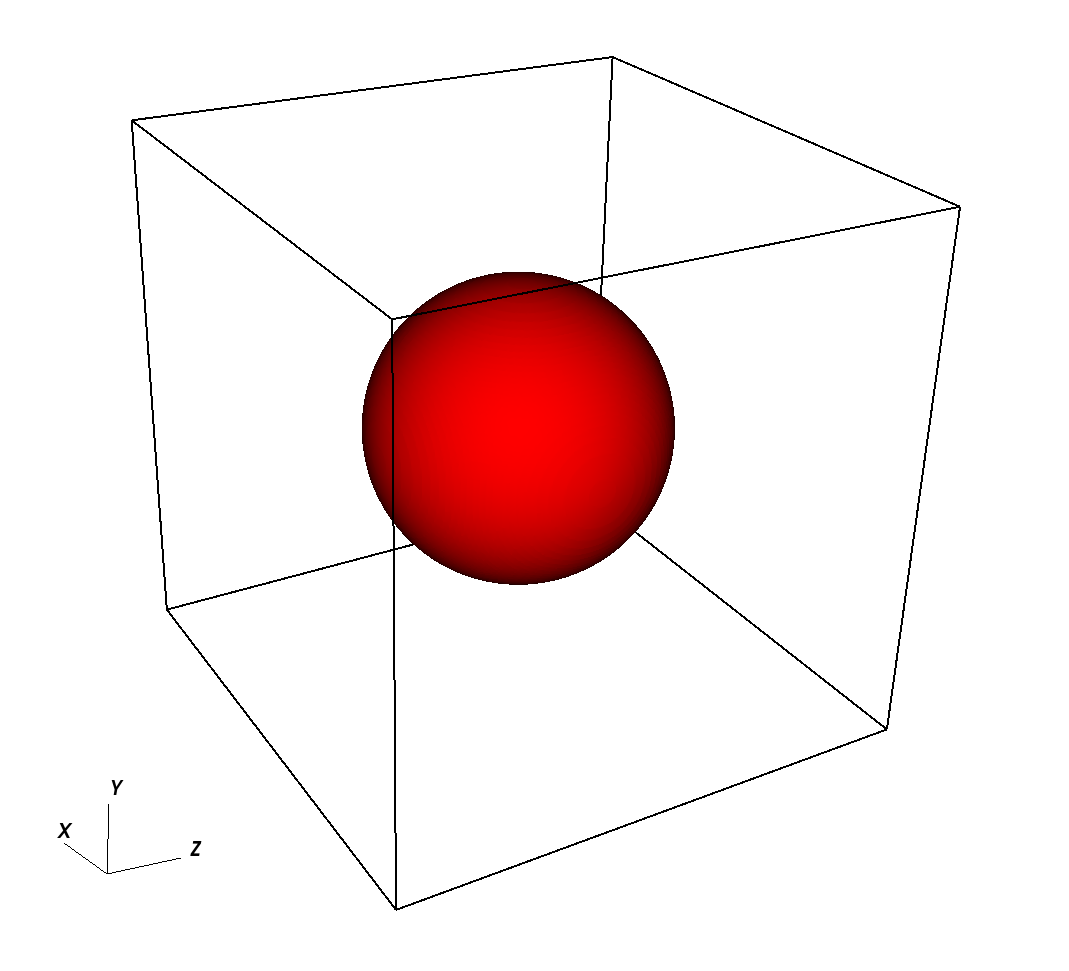}
\includegraphics[width=0.25\textwidth]{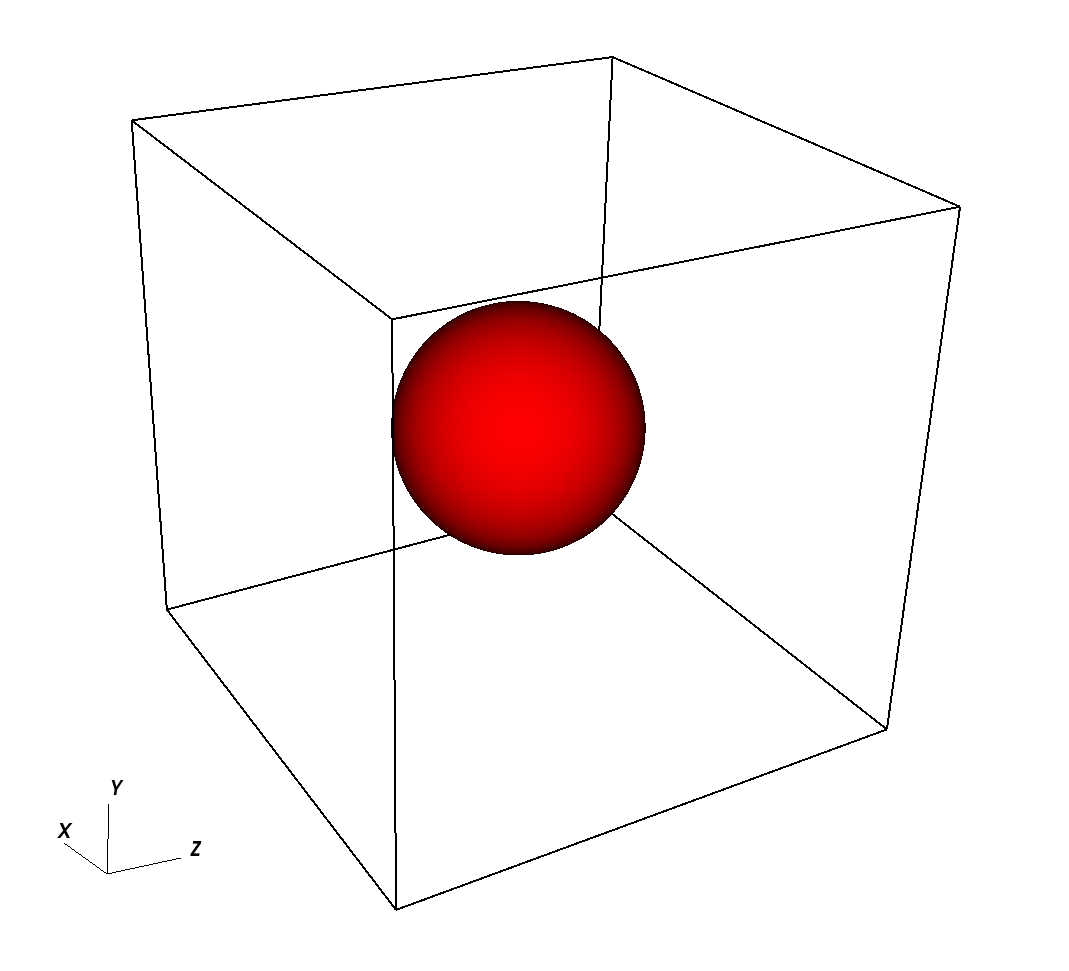}
\includegraphics[width=0.25\textwidth]{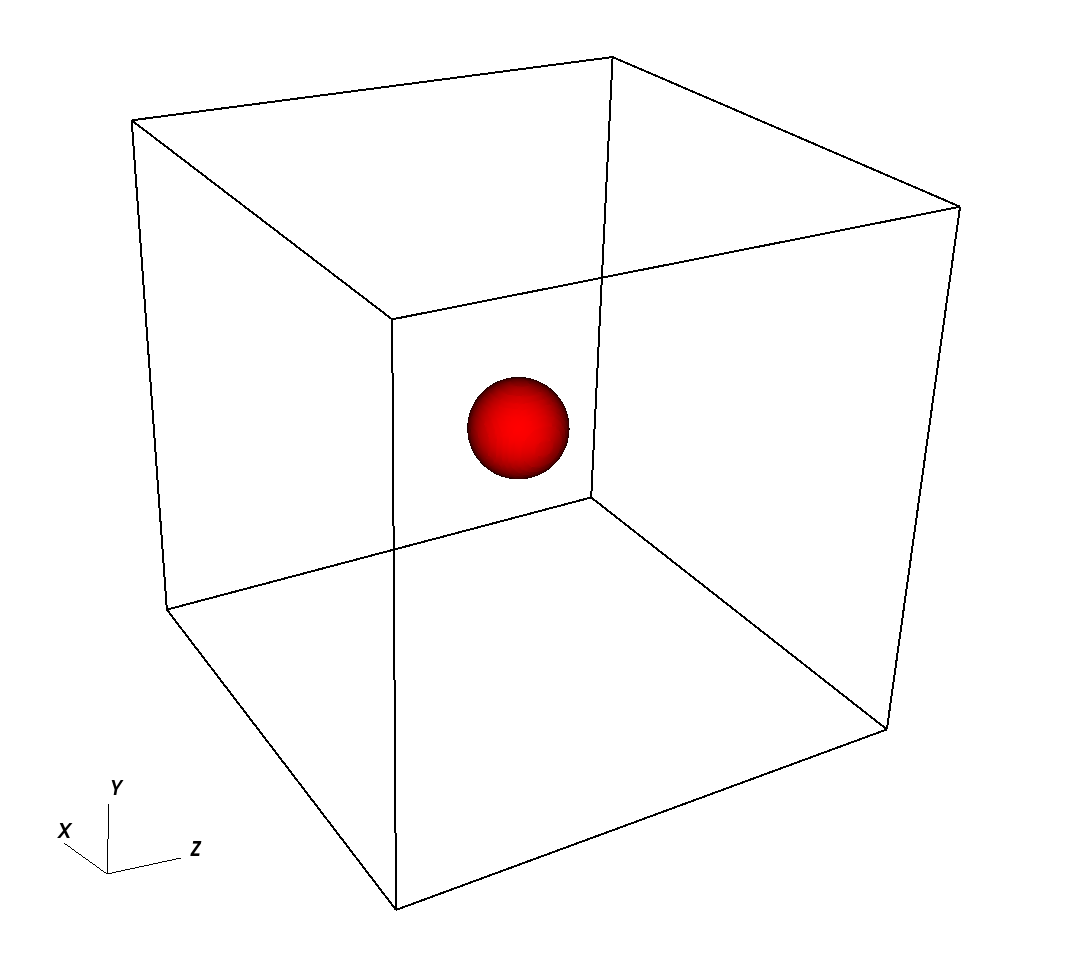}
}

\subfigure[Numerical errors ($u_{real}-u_{pred}$)  at $t=0, 0.2,0.5, 1$]{
\includegraphics[width=0.25\textwidth]{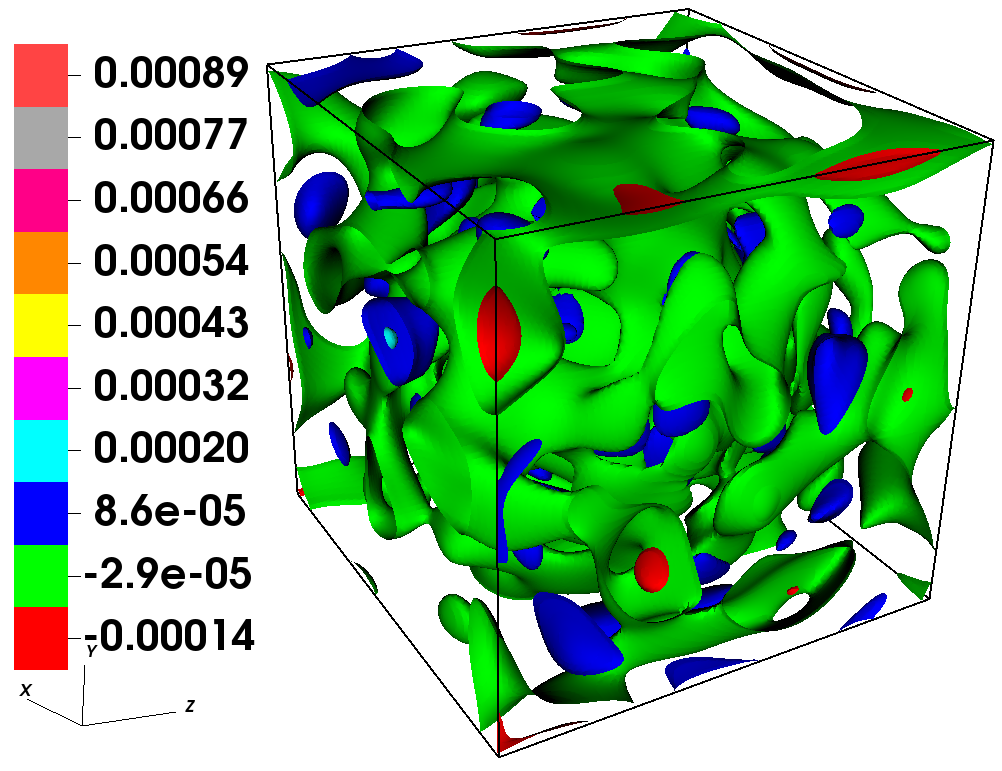}
\includegraphics[width=0.25\textwidth]{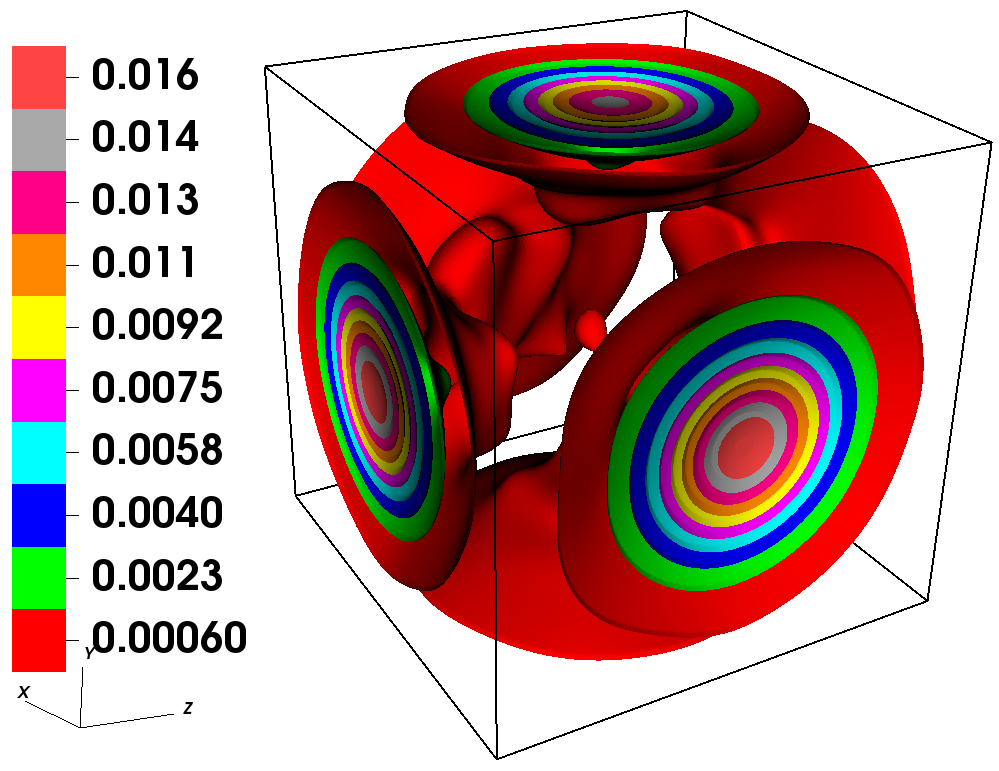}
\includegraphics[width=0.25\textwidth]{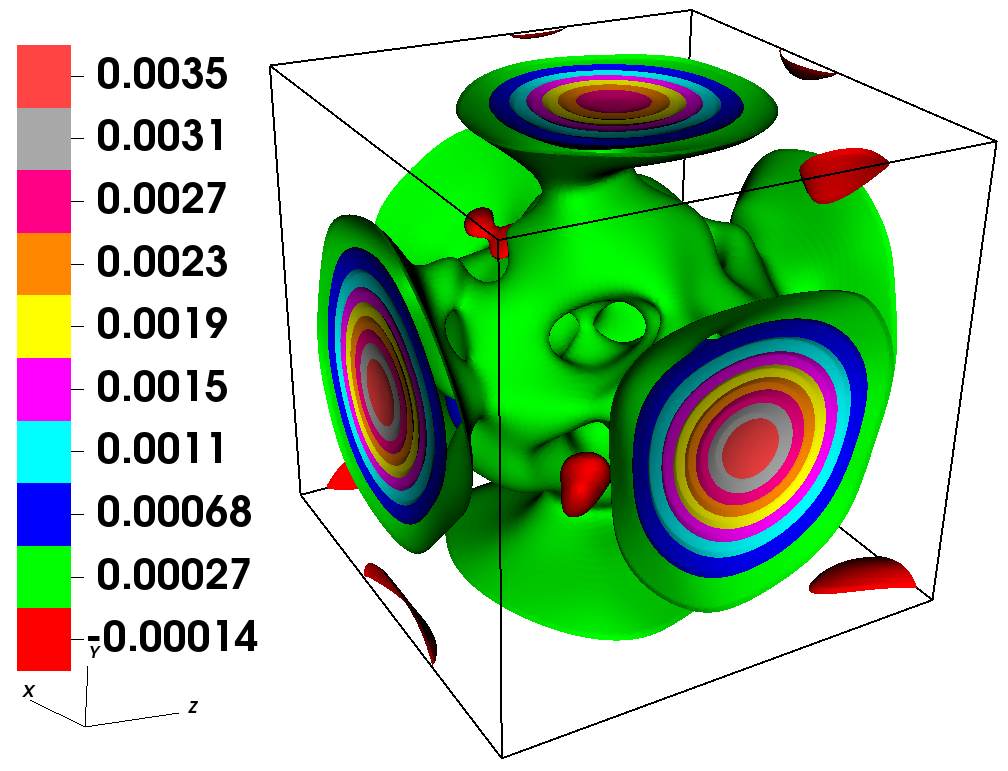}
\includegraphics[width=0.25\textwidth]{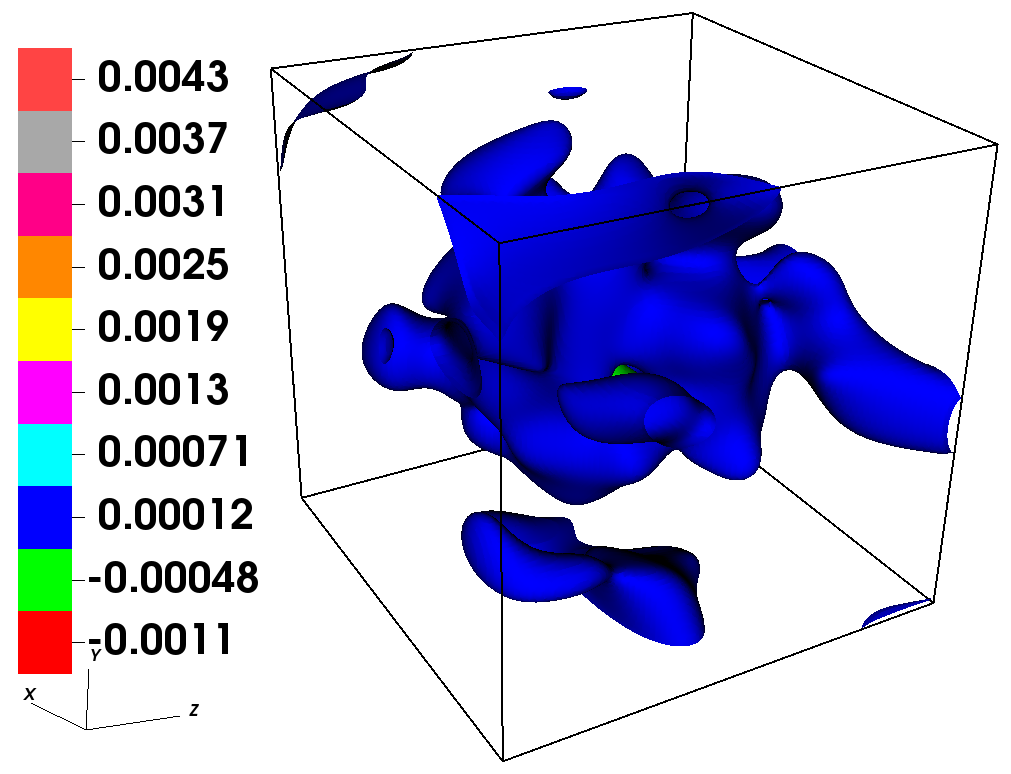}
}

\caption{Numerical approximations of 3D Allen-Cahn benchmark problem using the time-adaptive approach II. The predicted numerical solutions and their numerical errors at various time steps are shown. In this example, we chose $\Delta t=0.5$, and $C_0=10^3$   the neural network is $\cU: (x,y, z, t) \rightarrow \cU(x,y,z, t)$ with 6 hidden layers, 128 nodes per hidden layer. For each neural network, we used 100 epochs with batch size 32 for the Adam training, followed with a L-BFGS-B optimizer.}
\label{fig:AC3d}
\end{figure}

Similarly, the strategy can be used to solve the 3D benchmark problem as well. In this case, we chose $\Omega = [0, 1]^3$, with parameters $\lambda = 10$, $\varepsilon=0.05$. the initial profile for $\phi$ was chosen as
\beq
\phi(x, y, t= 0) = \tanh( \frac{ 0.35 - \sqrt{(x-0.5)^2 +(y-0.5)^2 + (z-0.5)^2}}{2 \varepsilon}).
\eeq
Then, the solutions for $t \in [0, 1]$ predicted by the time-adaptive approach II are summarized in Figure \ref{fig:AC3d}, where we observed that the neural network provided accurate predictions.

\begin{figure}[H]
\centering
\subfigure[Predicted numerical solutions ($u_{pred}$) at $t=0, 0.2, 0.5, 1$]{
\includegraphics[width=0.25\textwidth]{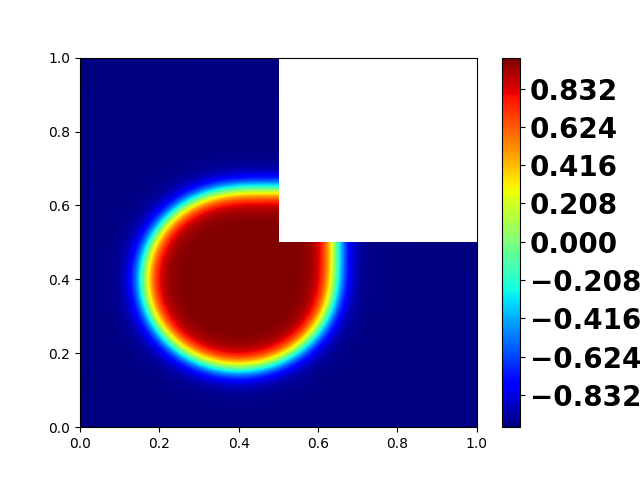}
\includegraphics[width=0.25\textwidth]{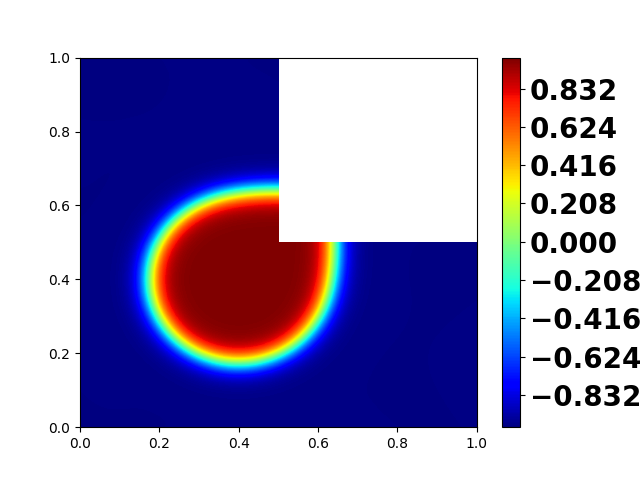}
\includegraphics[width=0.25\textwidth]{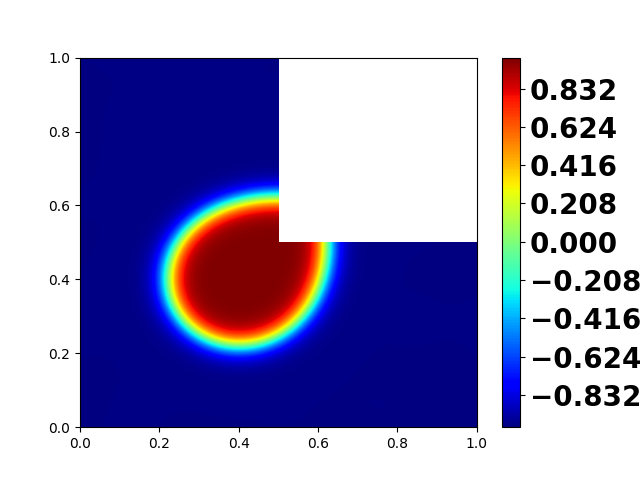}
\includegraphics[width=0.25\textwidth]{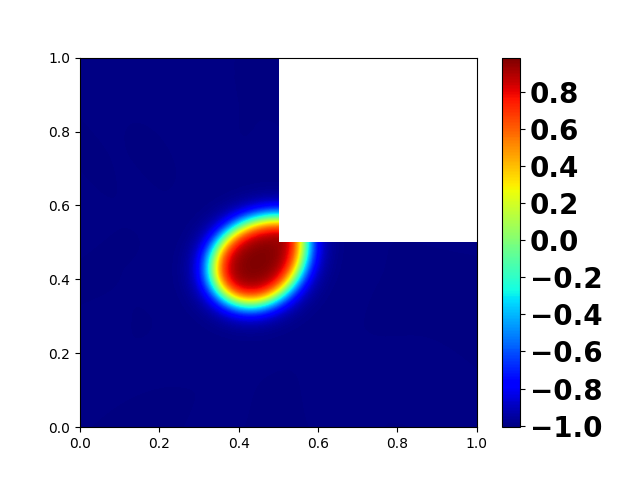}
}

\subfigure[Numerical errors ($u_{real}-u_{pred}$) at $t=0, 0.2, 0.5, 1$]{
\includegraphics[width=0.25\textwidth]{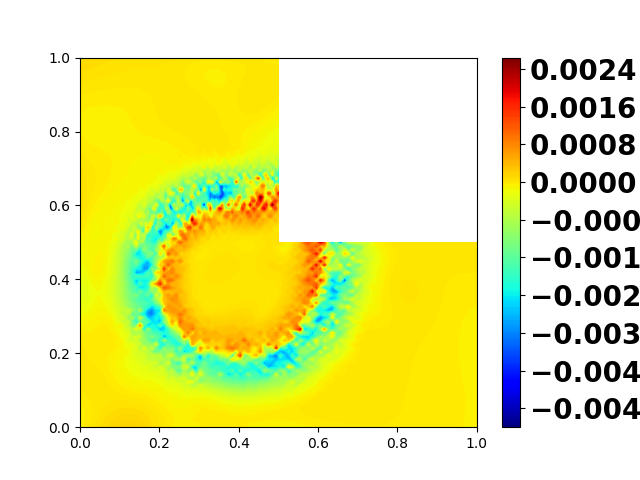}
\includegraphics[width=0.25\textwidth]{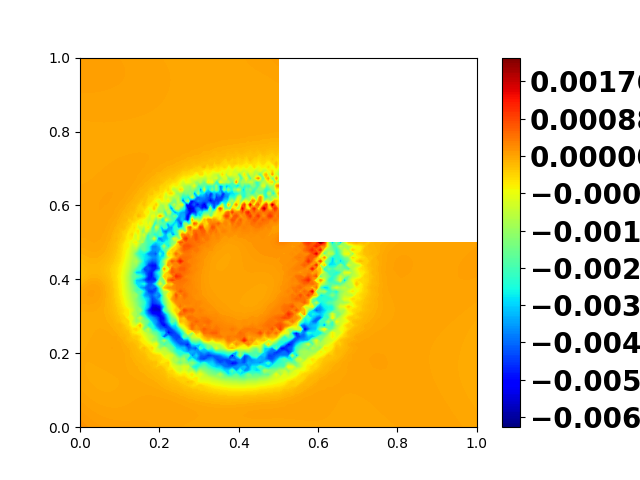}
\includegraphics[width=0.25\textwidth]{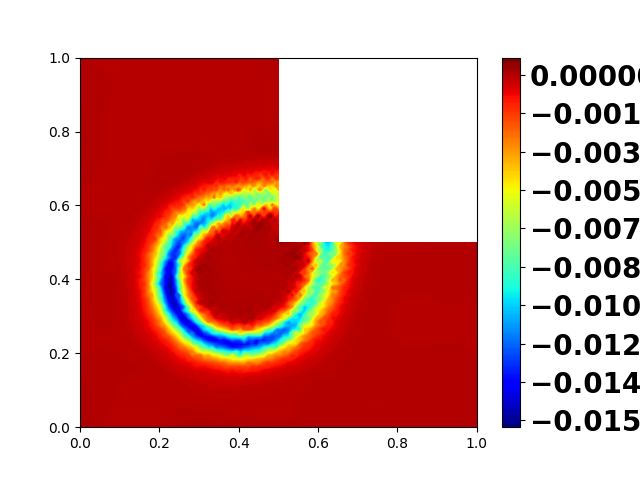}
\includegraphics[width=0.25\textwidth]{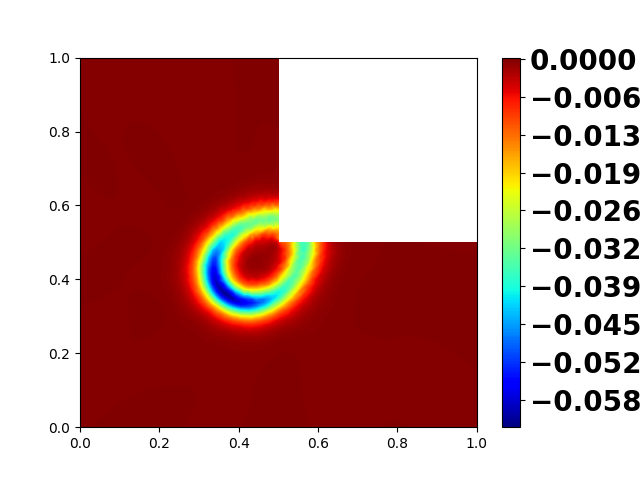}
}
\caption{Numerical approximations of the Allen-Cahn equation in a L-shape domain using the time-adaptive approach II. In this example, the neural network is $\cU: (x,y,z,t) \rightarrow \cU(x,y,z,t)$ with 6 hidden layers, 128 nodes per hidden layer. We chose time step $\Delta t = 0.1$ and $C_0=10^3$. For each neural network, we used 100 epochs with batch size 32 for the Adam training, followed with a L-BFGS-B optimizer.}
\label{fig:AC2d-Lshape}
\end{figure}

Then, to further demonstrate the effectiveness of the time-adaptive approach II, we solved the Allen-Cahn equation in a complex geometry, an L shape domain defined as $\Omega:= \{ 0 \leq x \leq 1 , 0 \leq y \leq 0.5\} \cup \{ 0 \leq x \leq 0.5 , 0.5 \leq y \leq 1\}$, with homogeneous Neumann boundary condition $\partial_\mathbf{n} \phi |_{\partial \Omega} = 0$. And we chose $\lambda=50$, $\varepsilon=0.025$, along with the initial condition
\beq
\phi(x, y, t=0) = \tanh \frac{ 0.25 - \sqrt{(x-0.4)^2+(y-0.4)^2}}{2\varepsilon}, \quad (x, y) \in \Omega.
\eeq 
The numerical results are shown in Figure \ref{fig:AC2d-Lshape}, where we observed the neural network predicted accurate dynamics of drop shrinking with relatively small errors.

\subsection{Solving the Cahn-Hilliard equation}
In the previous sub-section, we conducted a detailed study of the improved PINN on solving the Allen-Cahn equation. It indicated the proposed strategies have significantly improved the accuracy and convergence of the PINN.

Next we moved onto the Cahn-Hilliard Equation, which has higher order derivatives and is known to be harder to solve than the Allen-Cahn equation. Mainly, we focused on the following specific form
\begin{equation} \label{eq:CH-eq}
\begin{split}
&u_t - ( \gamma_2(u^3-u)-\gamma_1 u_{xx})_{xx}=0,\quad x \in [-1,1],\quad t\in [0,1], \\
&u(0,x) = -cos(2\pi x), \\
&u(t,-1)=u(t,1), \\
&u_x(t,-1)=u_x(t,1),
\end{split}
\end{equation}
with $\gamma_1, \gamma_2$ the model parameters.
For the following trials we used this equation with parameters fixed as $\gamma_2 = 0.01$, and $\gamma_1 =10^{-6}$.
We tested the same adaptive time method here as with Allen-Cahn with the regular adjustments to the $f$-network. 

Notice the fact the back-propagation with high order derivatives are extremely expensive to calculate. To overcome this computational deficiency, we introduced an intermediate $\mu$-network. Then, for solving the Cahn-Hilliard equation, we modified the networks as follows: the ($u$, $\mu$)-network as
\beq
\mathcal{U}: (x, t) \rightarrow [\mathcal{U}(x,t), \mu(x,t)],
\eeq 
and the $f$-network as
\beq
\mathcal{F}: (x, t) \rightarrow \cU_t(x,t) - \mu_{xx}.
\eeq
And the loss function was defined as
\beq 
MSE = MSE_u+MSE_b + MSE_f,
\eeq 
where the three terms were defined as 
\beq 
\begin{split}
&MSE_u=\frac{1}{N_u}\sum_{i=1}^{N_u}|\cU(0,x^i_u)-u^i|^2,  \\
& MSE_b=\frac{1}{N_b}\sum_{i=1}^{N_b}|\cU(t_b^i,x_u^i)-\cU(t_b^i,x_l^i)|^2+|\cU_x(t_b^i,x_u^i)-\cU_x(t_b^i,x_l^i)|^2,\\
&MSE_f=\frac{1}{N_f}\sum_{i=1}^{N_f}|\cF(t^i_f,x^i_f)|^2 + | \mu(t_f^i, x_f^i) -  \gamma_2(\cU(t_f^i, x_f^i)^3 - \cU(t_f^i, x_f^i)) + \gamma_1 \cU(t_f^i, x_f^i)_{xx}|^2.
\end{split}
\eeq 
This approach helps reducing the need to take the high order derivatives, and it turns out to significantly speed up computation and improve accuracy.

\begin{figure}
\centering
\includegraphics[scale=.65]{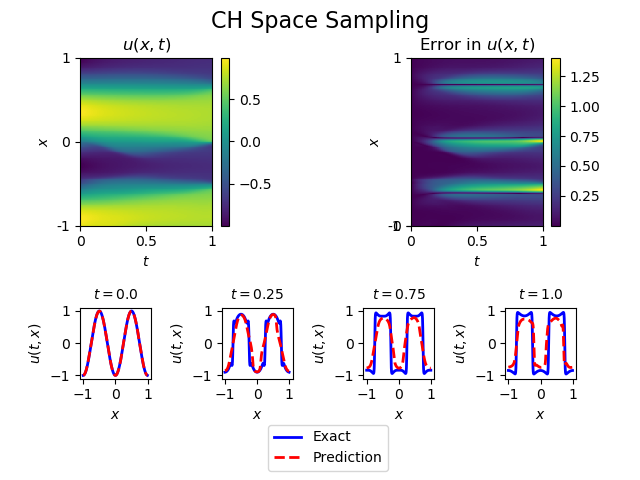}
\caption{Solutions of the Cahn-Hillard equation learned using the improved PINN with sampling collocation points alone. This figure indicates that training the PINN with sampling collocation points alone is not enough to achieve convergence for this problem even after many re-sampling iterations.}
\label{fig:CH-1}
\end{figure}

This Cahn-Hilliard equation in \eqref{eq:CH-eq} was tested using the approaches of sampling collocation points in Section \ref{sec:adaptive-sample} and time sampling strategies in Section \ref{sec:adaptive-in-time}. We observed that the trial of sampling collocation points alone does not converge, as the results illustrated in Figure \ref{fig:CH-1}. The neural network even failed to learn the solution well near $t=0$, and the error continued to accumulate throughout the entire domain. Like the trials for solving the AC equation, this trial used the same network architecture, weighted loss function, and the number of training data points. Each resampling iterations performed the Adam optimizer and the L-BFGS-B optimizer. The Adam optimizer learning rate is set to $0.001$, and the max iteration was set to $20,000$. Even with 20 resampling iterations, the solution did not converge. 


\begin{figure}[H]
\centering
\includegraphics[scale=.65]{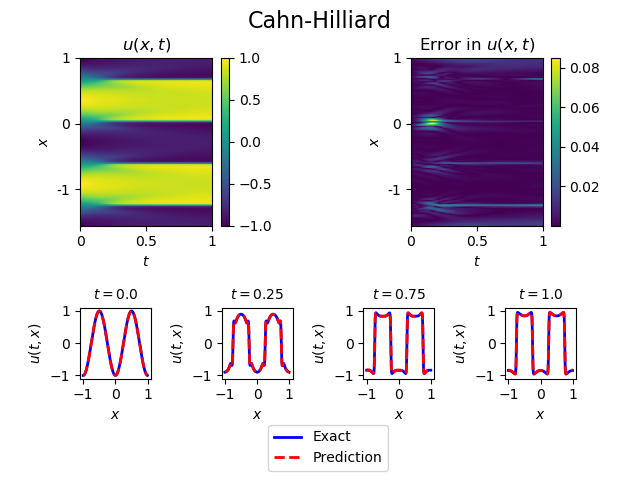}
\caption{Solutions of the Cahn-Hillard equation learned using the improved PINN with time-adaptive approach I in Section \ref{sec:adaptive-in-time} and adaptive re-sampling collocation points in Section \ref{sec:adaptive-sample}. 
The PINN with sampling collocation points alone was unable to learn the solution of the Cahn-Hilliard equation accurately. With the addition of adaptive time sampling, the improved PINN was able to learn an accurate solution with very small numerical errors.}
\label{fig:CH-2}
\end{figure}

In the meanwhile, when using the time sampling approaches in Section \ref{sec:adaptive-in-time}, we observed noticeable improvement. The neural network was able to learn the solution well near $t=0$, and by gradually allowing collocation points to be sampled at later times, the neural network maintained this accuracy across the whole domain.
Finally, we present the results for the trials run on solving the Cahn-Hilliard equation, which was the most difficult problem to solve thus far, in Figure \ref{fig:CH-2}. From our observation, only methods that involved adaptive sampling in space and time were able to learn the solutions accurately. In Figure \ref{fig:CH-2}, it shows the first trial using the time-adaptive approach I  on the Cahn-Hilliard equation, we obtained our best result with a relative $l_2$ error of $9.51e-3$.

\begin{figure}[H]
\centering
\subfigure[predicted numerical solutions ($u_{pred}$) at $t=0, 0.25, 0.5$ and $1$]{
\includegraphics[width=0.25\textwidth]{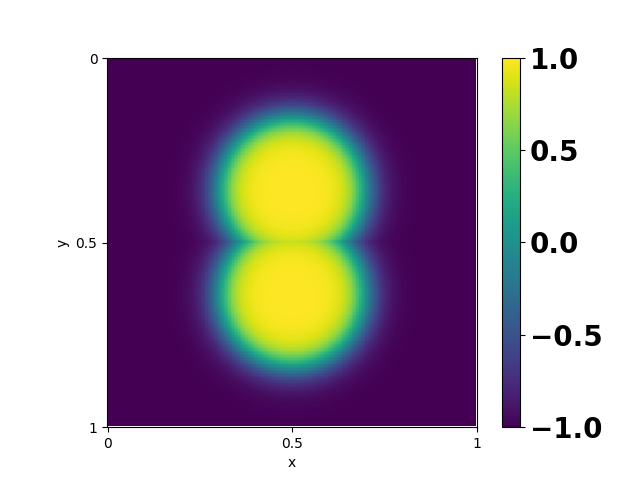}
\includegraphics[width=0.25\textwidth]{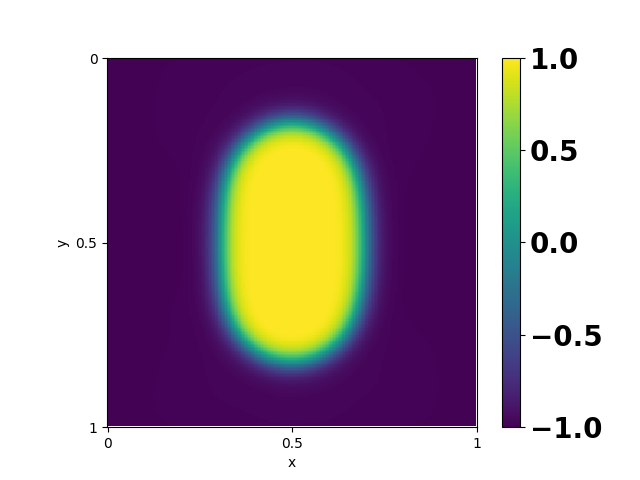}
\includegraphics[width=0.25\textwidth]{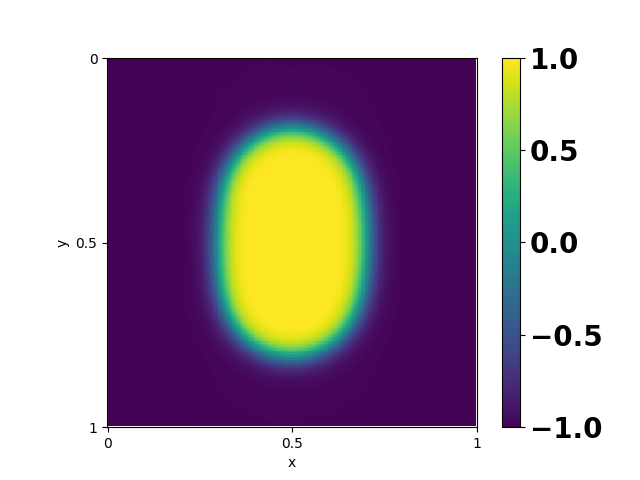}
\includegraphics[width=0.25\textwidth]{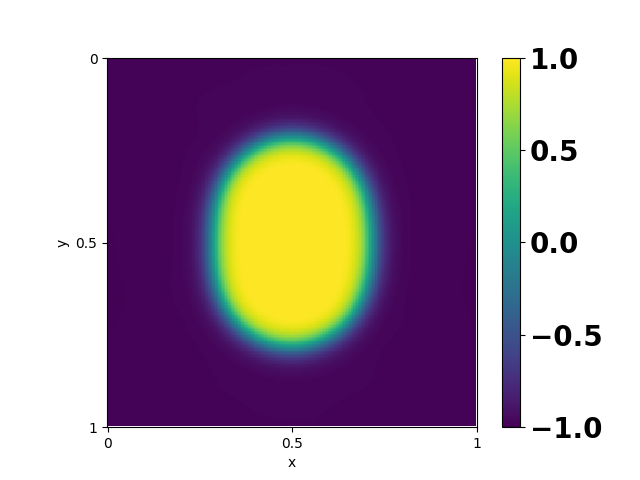}
}

\subfigure[Numerical errors ($u_{real}-u_{pred}$) at $t=0, 0.25, 0.5$ and $1$]{
\includegraphics[width=0.25\textwidth]{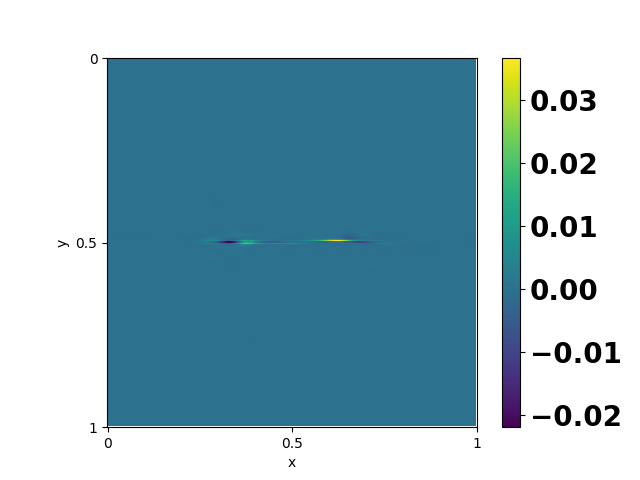}
\includegraphics[width=0.25\textwidth]{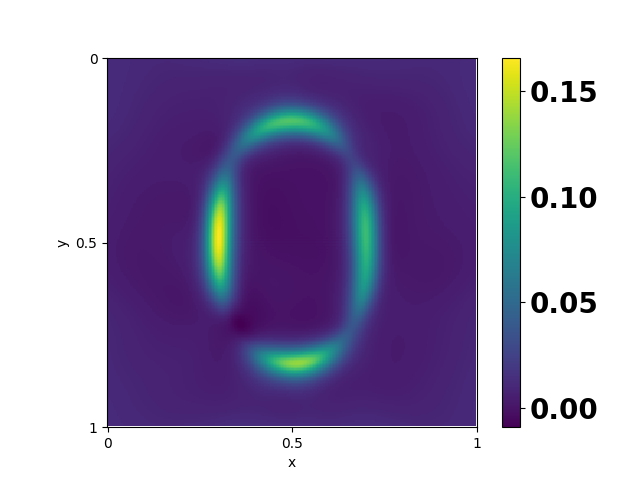}
\includegraphics[width=0.25\textwidth]{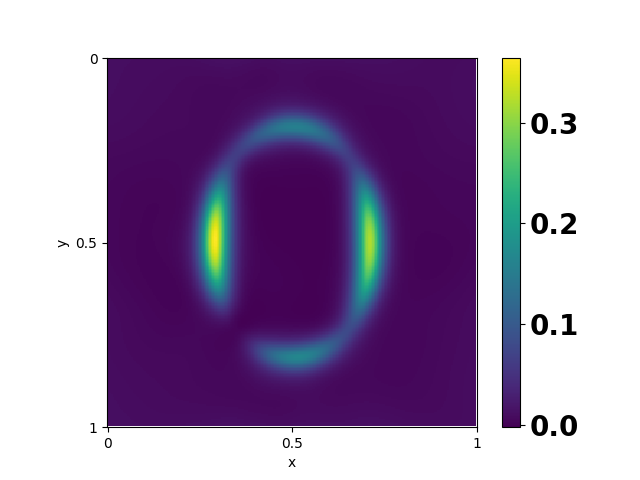}
\includegraphics[width=0.25\textwidth]{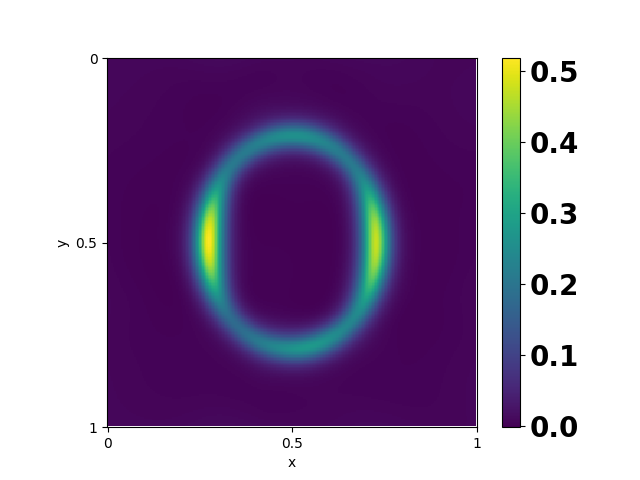}
}
\caption{Numerical approximations of the 2D Cahn-Hilliard equation using the time-adaptive approach II. In this example, the neural network is $\cU: (x,y,t) \rightarrow [\cU(x,y,t), \mu(x,y,t)]$ with 6 hidden layers, 128 nodes per hidden layer. We chose time step $\Delta t = 0.2$ and $C_0=10^3$. For each neural network, we used 100 epochs with batch size 32 for the Adam training, followed with a L-BFGS-B optimizer.}
\label{fig:CH2d}
\end{figure}

Next, we studied the benchmark problem for the Cahn-Hilliard model, which reads as
\beq
\partial_t \phi  = \lambda \Delta ( -\varepsilon^2 \Delta \phi + \phi^3 -\phi), \quad  (x, y) \in \Omega, t\geq 0.
\eeq
In this example, we chose the domain $\Omega:=[-1, 1]^2$, the parameters $\lambda=1$, $\varepsilon=0.05$, and the initial profile for $\phi$ as
\beq
\phi(x, y, t= 0) = \max (\tanh \frac{ r- R_1}{2\varepsilon}, \tanh \frac{r -R_2}{2\varepsilon}), \quad 
\eeq
where $r=0.4$, $R_1 = \sqrt{(x-0.7r)^2 + y^2}$, and $R_2 = \sqrt{(x+0.7r)^2 + y^2}$. This problem was solved for $t\in[0, 1]$ by the time adaptive approach II, and the predicted solutions are summarized in Figure \ref{fig:CH2d}(a) and the numerical errors (the difference between the real solutions and the predicted solutions) are summarized in Figure \ref{fig:CH2d}(b). We observed the numerical solutions predicted by the adaptive approach II can capture the the bubble merging accurately. 

\section{Conclusion}
In this paper, we have introduced several strategies to improve the approximating capability of the physics informed neural network (PINN). Then we used the improved PINN to solve the phase field equations of increased complexity. Even though we focused on the problem of solving the Allen-Cahn equation and the Cahn-Hilliard equation, the improved PINN could readily be used to solve other difficult phase field equations as well. 

Space sampling opened the door to other ideas of adaptive sampling. Space sampling uses the $f$-network predictions to pinpoint areas to focus data points on. Time sampling uses knowledge of differential equations to chose areas to focus on. We saw how both the value of the loss function and the values of the $f$-network predictions could help in determining if a network has learned a solution well. Using this information can help the network in making a decision, such as whether to use more collocation points or whether to focus on a smaller time domain.
We also saw merit in both of the time sampling methods. The time-adaptive approach I proved better than just space sampling alone. The time-adaptive approach II took a different angle that has the potential for even higher accuracy and learning the solution faster. It uses individual networks that can focus on a smaller problem domain. A potential downside with more difficult equations is that if the solution gets off on a time interval, all the intervals after that will propagate that error. It is important to learn the solution well on a time interval before moving on to the next one. This method may also be helpful when working on problems with larger time domains.

We have seen how adapting classical mathematical techniques and principles can help us find useful approaches for designing and training artificial deep neural networks. Roughly speaking, the best performance is obtained by using a combination of all of the techniques presented. More simple methods such as mini-batching, and adding weights in loss function are useful, especially when they are combined with more powerful adaptive sampling methods.
This research has focused mainly on the problem of {\it solution} of differential equations. The next step of our future work is to test these approaches on the {\it discovery} of differential equations. The same neural network architectures can be used with the addition of a few more learnable/trainable parameters. 

\section*{Acknowledgments}
The authors would like to acknowledge the support from NSF-DMS-1816783 
and NVIDIA Corporation for their donation of a Quadro P6000 GPU for conducting some of the numerical simulations in this paper.


\end{document}